\documentclass{article}
\usepackage{amsthm, amsmath,amssymb, amscd}
\usepackage{txfonts}
\usepackage[dvips]{graphicx, color}
\usepackage{latexsym,amsxtra,amsfonts,verbatim}
\usepackage{lipsum}     
\usepackage[makeroom]{cancel}

\usepackage[dvips]{graphicx, color}
\usepackage[title]{appendix}

\newcommand{\xdownarrow}[1]{%
  {\left\downarrow\vbox to #1{}\right.\kern-\nulldelimiterspace}
}
\usepackage[new]{old-arrows}

\usepackage{bbm}

\theoremstyle{plain}
\newtheorem{thm}{Theorem}[section]
\newtheorem{lem}[thm]{Lemma}

\newtheorem{prop}[thm]{Proposition}
\newtheorem{cor}[thm]{Corollary}

\theoremstyle{definition}
\newtheorem{defn}{Definition}[section]
\newtheorem{conj}{Conjecture}[section]
\newtheorem{exam}{Example}[section]

\theoremstyle{remark}
\newtheorem*{rmk}{Remark}

\newcommand{\cE}{{\mathcal E}}
\newcommand{\cF}{{\mathcal F}}

\newcommand{\cL}{{\mathcal L}}

\newcommand{\C}{{\Bbb C}}

\newcommand{\F}{{\Bbb F}}

\renewcommand{\O}{{\Bbb O}}

\newcommand{\R}{{\Bbb R}}

\newcommand{\Z}{{\Bbb Z}}

\def\O{\mathcal O}

\def\lra{\longrightarrow}

\def\={\:=\:}  \def\+{\,+\,}

\def\a{\alpha} \def\b{{\beta}}  \def\ba{\overline\a}

\newcommand{\frZ}{{\frak  Z}}

\def\be{\begin{equation}}   \def\ee{\end{equation}}
\def\bes{\begin{equation*}}   \def\ees{\end{equation*}}
\def\ba{\begin{aligned}}   \def\ea{\end{aligned}}
\def\bc{\begin{cases}}   \def\ec{\end{cases}}
\def\bp{\begin{proof}}   \def\ep{\end{proof}}

\newcommand{\Res}{\mathrm{Res}}

\newcommand{\Aut}{\mathrm{Aut}}

\def\SL{\mathrm{SL}}

\newcommand{\ov}{\overline}
\def\qqan{\qquad\mathrm{and}\qquad}
\def\qan{\quad\mathrm{and}\quad}

\newcommand{\dis}{\displaystyle}

\def\smm{\smallsetminus}

\def\Ga{\Gamma}

\def\SL{\mathrm{SL}}

\def\Llra{\Longleftrightarrow}
\def\lra{\longrightarrow}

\def\ov{\overline}
\def\De{\Delta}

\def\llra{\longleftrightarrow}

\def\bbm1{\mathbbm 1}

\def\De{\Delta}

\def\wh{\widehat}
\def\be{\begin{equation}}   \def\ee{\end{equation}}
\def\bes{\begin{equation*}}   \def\ees{\end{equation*}}
\def\bea{\begin{equation}\begin{aligned}}   
\def\eea{\end{aligned}\end{equation}}

\def\wt{\widetilde}

\def\Pic{\mathrm{Pic}}

\def\bm{\begin{matrix}}
\def\em{\end{matrix}}
\def\bpm{\begin{pmatrix}}
\def\epm{\end{pmatrix}}

\def\Z{\mathbb Z} \def\R{\mathbb R}    
\def\={\;=\;}  \def\+{\,+\,}  \def\C{\Bbb C}      
\def\Z{\Bbb Z}  \def\F{\Bbb F}

\begin{document}

\title{\bf  Derived Zeta Functions\\
 for Curves over Finite Fields} 
\author{Lin WENG}  
\date{}
\maketitle
\tableofcontents
\begin{abstract}
For each $(m+1)$-tuple ${\bf n}_m=(n_0,n_1,\ldots,n_m)$ of positive integers, the ${\bf n}_m$-derived zeta function $\wh\zeta_{X,\F_q}^{\,({\bf n}_m)}(s)$ is defined for a curve $X$ over $\F_q$, motivated by the theory of rank $n$ non-abelian zeta function $\wh\zeta_{X,\F_q;n}(s)$ of $X/\F_q$. This derived zeta function satisfies standard zeta properties such as the rationality, the functional equation and admits only two singularities, namely, two simple poles at $s=0,1$, whose residues are given by 
the ${\bf n}_m$-derived beta invariant $\b_{X,\F_q}^{\,({\bf n}_m)}$ for which the Harder-Narasimhan-Ramanan-Desale-Zagier type formula holds. In particular, similar to the Artin Zeta function of $X/\F_q$, this 
${\bf n}_m$-derived Zeta function for  $X$ over $\F_q$ is a ratio of a degree $2g$ polynomial $P_{X,\F_q}^{({\bf n}_m)}$ in $T_{{\bf n}_m}=q^{-s\prod_{k=0}^mn_k}$ by $(1-T_{{\bf n}_m})(1-q_{{\bf n}_m}T_{{\bf n}_m})T_{{\bf n}_m}^{g-1}$ with $q_{{\bf n}_m}=q^{\prod_{k=0}^mn_k}$. Indeed, $\wh \zeta_{X,\F_q}^{({\bf n}_{m})}(s)=\wh Z_{X,\F_q}^{({\bf n}_{m})}(T_{{\bf n}_{m}})$ are given by
$$\ba
&
\left(\sum_{\ell=0}^{g-2}\a_{X,\F_q}^{({\bf n}_{m})}(\ell)\Big(T_{{\bf n}_{m}}^{\ell-(g-1)}+q_{{\bf n}_{m}}^{(g-1)-\ell}T_{{\bf n}_{m}}^{(g-1)-\ell}\Big)
+\a_{X,\F_q}^{({\bf n}_{m})}(g-1))\Big)\right)+\frac{(q_{{\bf n}_{m}}-1)T_{{\bf n}_{m}}\b_{X,\F_q}^{({\bf n}_{m})}}{(1-T_{{\bf n}_{m}})(1-q_{{\bf n}_{m}}T_{{\bf n}_{m}})}\\
\ea$$
for some ${\bf n}_m$-derived alpha invariants $\Big\{\a_{X,\F_q}^{({\bf n}_{m})}(\ell)\Big\}_{\ell=0}^{g-1}$ of $X/\F_q$. Furthermore, when $X$ restrict to an elliptic curve, or when ${\bf n}_m=(2,2,\ldots 2)$, established is the ${\bf n}_m$-derived Riemann hypothesis claiming that all zeros of $\wh \zeta_{X,\F_q}^{({\bf n}_{m})}(s)$ lie on the central line $\Re(s)=\frac{1}{2}$. In addition,  formulated is the Positivity Conjecture claiming that the above  ${\bf n}_m$-derived alpha and beta invariants  are all strict positivity. This Positivity Conjecture is the key to control our ${\bf n}_m$-derived zetas.
\end{abstract}

\section{${\bf n}_m$-Derived Zeta Functions for Curves over Finite Fields}
In this section, we define inductively the ${\bf n}_m$-derived zeta functions for curves over finite fields associated to $(m+1)$-tuples ${\bf n}_m=(n_0,\ldots,n_m)$ of positive integers $n_0,n_1,\ldots,n_m$.

\subsection{Rank $n$ zeta function}
In this subsections, as an initial step in the inductive process to introduce the ${\bf n}_m$-derived zeta functions for curves over finite fields, we recall some basic structures of rank $n$ non-abelian zeta functions for these curves.

 Let $X$ be an integral  regular projective curve of genus $g$ over a finite field $\F_q$ with $q$ elements. We define the {\it rank $n$ non-abelian zeta function $\wh\zeta_{X,\F_q;n}(s)$  of $X/\F_q$ } by
\be
\wh\zeta_{X,\F_q;n}(s):=\sum_{m\geq 0}\sum_{\cE}\frac{q^{h^0(X,\cE)}-1}{\#\Aut(\cE)}(q^{-s})^{\chi(X,\cE)}\qquad\Re(s)>1.
\ee
Here $\cE$ runs through rank $n$ semi-stable $\F_q$-rational vector bundles over $X/\F_q$ of degree $mn$.
Tautologically, by applying the Riemann-Roch theorem, the cohomological duality and the vanishing theorem for semi-stable bundles over $X$, we have, for $Q=q^n, T=q^{-ns}$,
$$\ba
&\wh Z_{X,F_q;n}(T):=\wh\zeta_{X,\F_q;n}(s)\\
=&
\sum_{m=0}^{g-2}\a_{X,\F_q;n}(nm)\Big(T^{m-(g-1)}+(QT)^{(g-1)-m}\Big)
+\a_{X,\F_q;n}(n(g-1))+\frac{(Q-1)T}{(1-T)(1-QT)}\b_{X,\F_q;n}(0)\\
=:&\frac{P_{X,\F_q;n}(T)}{(1-T)(1-QT)T^{g-1}}
\ea
$$
where \be
\a_{X,\F_q;n}(d)=\sum_{\cE}\frac{q^{h^0(X,\cE)}-1}{\#\Aut(\cE)}\qqan
\b_{X,\F_q;n}(d):=\sum_{\cE}\frac{1}{\#\Aut(\cE)}
\ee
where $\cE$ runs over all rank $n$ semi-stable $\F_q$-rational vector bundle of degree $d$ on $X$
and
$$\ba
&P_{X,\F_q;n}(T)\\
:=&\left(\sum_{m=0}^{g-2}\a_{X,\F_q;n}(nm)\Big(T^{m}+Q^{(g-1)-m}T^{2(g-1)-m}\Big)
+\a_{X,\F_q;n}\big(n(g-1)\big)T^{g-1}\right)(1-T)(1-QT)\\
&\hskip 8.60cm+(Q-1)T^g\b_{X,\F_q;n}(0)
\ea
$$
is a degree $2g$ polynomial in $T$ with rational coefficients.

\begin{exam} When $n=1$, we recover the (complete) Artin zeta function for $X$ over $\F_q$:
$$
\ba
\wh\zeta_{X,\F_q;1}(s)=&\sum_{d\geq 0}\sum_{\cL\in\Pic^d(X)}\frac{q^{h^0(X,\cL)}-1}{q-1}(q^{-s})^{d(\cL)-(g-1)}\\
=&q^{s(g-1)}\sum_{m\geq 0}\sum_{D\geq 0}N(D)^{-s}=:q^{s(g-1)}\zeta_{X,\F_q}(s)
=:\wh\zeta_{X,\F_q}(s)
\ea
$$
where $\Pic^d(X)$ denotes the degree $d$ Picard group of $X/\F_q$, and $D$ runs over effective divisors of degree $m$ on $X$.
\end{exam}
Furthermore, by the discussion, we have the following:
\begin{thm} [Zeta Facts \cite{WZeta}]\label{thm1.1} Let $X$ be an integral  regular projective curve of genus $g$ over $\F_q$. Then we have
\begin{enumerate}
\item [(1)] {\rm(Rationality)}
$$
\wh Z_{X,\F_q;n}(T):=\wh\zeta_{X,\F_q;n}(T)=\frac{P_{X,\F_q;n}(T)}{(1-T)(1-QT)T^{g-1}}
$$
is a rational function in $T$, and $P_{X,\F_q;n}(T)$ is a degree $2g$ polynomial with rational coefficients in $T$.
\item[(2)] {\rm (Functional Equation)}
$$
\wh\zeta_{X,\F_q;n}(1-s)=\wh\zeta_{X,\F_q;n}(s),
\quad{\rm
or\ equivalently,}\quad
\wh Z_{X,\F_q;n}(1/(QT))=\wh Z_{X,\F_q;n}(T)
$$
\item[(2)] {\rm (Singularities)}
The rank $n$  Zeta function
$Z_{X,\F_q;n}(T):=T^{g-1}\wh Z_{X,\F_q;n}(T)$ admits only two singularities, namely two simple poles at $T=1,1/Q$ whose residues are given by
$$
\Res_{T=1}\wh Z_{X,\F_q;n}(T)=-\Res_{T=1/Q}\wh Z_{X,\F_q;n}(T)=\beta_{X,\F_q;n}(0).
$$
\end{enumerate}
\end{thm}

Accordingly, by using the functional equation, we may write
\be
P_{X,\F_q;n}(T)=\a_{X,\F_q;n}(0)\prod_{i=1}^g(\a_{X,\F_q;i}T-1)(\b_{X,\F_q;i}T-1)=\a_{X,\F_q;n}(0)\prod_{i=1}^g(QT^2-a_{X,\F_q;n;i}T+1)
\ee
where $(\a_{X,\F_q;n;i},\b_{X,\F_qn;i})\ (1\leq i\leq g)$  denotes the naturally paired reciprocal roots of $P_{X,\F_q;n}(T)$ characterized by the following conditions
$$
\a_{X,\F_q;n;i}\cdot \b_{X,\F_q;n;i}=Q\qqan \a_{X,\F_q;i}+\b_{X,\F_q;i}=a_{X,\F_q;n;i}.
$$
For our own convenience, we often write
$$
\a_{X,\F_q;2g-i}:=\b_{X,\F_q;i}\qquad (i=1,\ldots, g).
$$
and write the polynomial $P_{X,\F_q;n}(T)$ as
\be
\ba
&P_{X,\F_q;n}(T)\\
=&\Biggl(\Bigg(\sum_{m=0}^{g-2}\a_{X,\F_q;n}'(nm)\Big(T^{m}+Q^{(g-1)-m}T^{2(g-1)-m}\Big)+\a_{X,\F_q;n}'(n(g-1))T^{g-1}\Bigg)(1-T)(1-QT)\\
&\hskip 6.80cm+(Q-1)T^g\b_{X,\F_q;n}'(0)\Biggr)\cdot \a_{X,\F_q;n}(0)\\
\ea
\ee
Here, we have set
\be
\a_{X,\F_q;n}'(d)=\frac{\a_{X,\F_q;n}(d)}{\a_{X,\F_q;n}(0)}\qqan \b_{X,\F_q;n}'(d)=\frac{\b_{X,\F_q;n}(d)}{\a_{X,\F_q;n}(0)}\qquad (\forall d\in \Z)
\ee
using the fact that 
$$
\a_{X,\F_q;n}(0)>\frac{q^{h^0(X,\O_X^{\oplus n})}-1}{\#\Aut(\O_X^{\oplus n})}>0,
$$
Obviously, the leading coefficient  and constant term of $P_{X,\F_q;n}(T)$ are given by
$\a_{X,\F_q;n}(0)Q^g$ and $\a_{X,\F_q;n}(0)$, respectively.

There is one more important zeta properties, at least conjecturally, for these zeta functions. Namely,
\begin{conj}[Rank $n$ Riemann Hypothesis] Let $X$ be an integral  regular projective curve of genus $g$ over $\F_q$. Then all zeros of the rank $n$ non-abelian zeta function $\wh\zeta_{X,\F_q;n}(s)$  of $X/\F_q$ lie on the central line $\Re(s)=\frac{1}{2}$.
\end{conj}
we end this subsection with the following rather easy elementary

\begin{lem}\label{lem1.2} Let $X$ be an integral  regular projective curve of genus $g$ over $\F_q$. Then the following conditions are equivalent:
\begin{enumerate}
\item [(1)] The rank $n$ Riemann Hypothesis holds for $\wh\zeta_{X,\F_q;n}(s)$.
\item [(2)] For all $i=1,\ldots, 2g$, 
$$
|\a_{X,\F_q;i}|=\sqrt Q.
$$
\item [(3)] For all $i=1,\ldots, g$,
$$
\ov{\a_{X,\F_q;n;i}}=\b_{X,\F_q;n;i}.
$$
\item[(4)] For all $i=1,\ldots g$,
$$
a_{X,\F_q;n;i}\in \R\qqan a_{X,\F_q;n;i}\in (-2\sqrt Q, 2\sqrt Q)
$$
\item[(5)]  For all $i=1,\ldots g$, 
$$
\a_{X,\F_q;n;i}\in \C\smm\R\qqan \b_{X,\F_q;n;i}\in \C\smm \R.
$$
\end{enumerate}
\end{lem}

\subsection{Special uniformity of zetas}
There is another new type of zeta function $\wh\zeta_{X,\F_q}^{\SL_n}(s)$ for curves $X$ over finite fields $\F_q$. These zeta functions are defined using the Lie structures of $\SL_n$, or better of the pair $(\SL_n,P_{n-1,1})$, where $P_{n-1,1}$ denotes the maximal parabolic subgroup of $\SL_n$, associated to the ordered partition $n=(n-1)+1$, consisting of matrices whose final row vanishes except for its  last entry, and the complete Artin zeta function of $X/\F_q$. The so-called special uniformity of zeta functions for curves over finite fields claims that the geometrically defined rank $n$-non-abelian zeta function
$\wh\zeta_{X,\F_q;n}(s)$ coincides with the Lie theoretically defined $\SL_n$ zeta function 
$\wh\zeta_{X,\F_q}^{\SL_n}(s)$. This special uniformity of zetas was conjectures in \cite{W3} and proved in \cite{WZ2}, with the help of \cite{MR}.

 \begin{thm}[Special Uniformity of Zetas. Theorem 1 of \cite{WZ2}]\label{thm1.3} For an integral  regular projective curve $X$ of genus $g$
 over $\F_q$, we have, for $n\geq 2$,
 $$
 \ba
 \wh\zeta_{X,\F_q;n}(s)=&\wh\zeta_{X,\F_q}^{\SL_n}(s)=q^{\binom{n}{2}(g-1)}\sum_{a=1}^{n}\sum_{\substack{k_1,\ldots,k_p>0\\ k_1+\ldots+k_p=n-a}}\frac{\wh v_{X,\F_q,k_1}\ldots\wh v_{X,\F_q,k_p}}{\prod_{j=1}^{p-1}(1-q^{k_j+k_{j+1}})} \frac{1}{(1-q^{ns-n+a+k_{p}})}\\
 &\hskip 1.80cm\times\wh\zeta_{X,\F_q}(ns-n+a)\sum_{\substack{l_1,\ldots,l_r>0\\ l_1+\ldots+l_r=a-1}}
 \frac{1}{(1-q^{-ns+n-a+1+l_{1}})}\frac{\wh v_{X,\F_q,l_1}\ldots\wh v_{X,\F_q,l_r}}{\prod_{j=1}^{r-1}(1-q^{l_j+l_{j+1}})}\\
 \ea
 $$
 where $\dis{\wh\nu_{X,\F_q;n}:=\prod_{k=1}^n\wh\zeta_{X,\F_q}(k)}$ with $\,\wh\zeta_{X,\F_q}(1):=\Res_{T=1}\wh Z_{X,\F_q}(T)$.
 \end{thm}
 
 The importance of this special uniformity of zetas should never be underestimated. For examples, it has been used in \cite{WRH} to establish the rank three Riemann hypothesis, and as to be seen below plays an important role in defining what we call ${\bf n}_m$-derived zeta functions for curves over finite fields.

\subsection{Definition of ${\bf n}_m$-derived zeta functions}\label{sec1.3}
 \begin{defn}\label{defn1.1}
Let $X$ be  an integral  regular projective curve  of genus $g$ over $\F_q$.  Fix an $(m+1)$-tuple ${\bf n}_m$ of (strictly) positive integers ${\bf n}_m=(n_0,n_1,\ldots,n_m)$. The ${\bf n}_m$-derived zeta functions, or simply, the level $m$-derived zeta function, for $X$ over $\F_q$  is defined inductively by

\noindent
(1) When $m=0$, for ${\bf n}_0=(n)$, we set
$$
\wh\zeta_{X,\F_q}^{\,({\bf n}_0)}(s):=\wh\zeta_{X,\F_q}^{\SL_n}(s);
$$

\noindent
(2) For all $m\geq 1$,  inductively, assume that the ${\bf n}_k$-derived zeta functions
$\wh\zeta_{X,\F_q}^{\,({\bf n}_k)}(s)$ have been defined for all $k\leq m-1$. Then, for
${\bf n}_m:=(n_0,n_1,\ldots,n_m)$
we set
$$
q_{{\bf n}_m}=q_{{\bf n}_{m-1}}^{n_m}\qqan T_{{\bf n}_m}=T_{{\bf n}_{m-1}}^{n_m}
$$
and
$$
\ba
\wh\zeta_{X,\F_q}^{\,({\bf n}_m)}(s):=&
q_{{\bf n}_{m-1}}^{\binom{n^{~}_m}{2}(g-1)}\sum_{a=1}^{n_m}\sum_{\substack{k_1,\ldots,k_p>0\\ k_1+\ldots+k_p=n_m-a}}\frac{\wh v_{X,\F_q,k_1}^{({\bf n}_{m-1})}\ldots\wh v_{X,\F_q,k_p}^{({\bf n}_{m-1})}}{\prod_{j=1}^{p-1}(1-q_{{\bf n}_{m-1}}^{k_j+k_{j+1}})} \frac{1}{(1-q_{{\bf n}_{m-1}}^{n_ms-n_m+a+k_{p}})}\\
 &\hskip 0.20cm\times\wh\zeta_{X,\F_q}^{\,({\bf n}_{m-1})}(n_ms-n_m+a)\sum_{\substack{l_1,\ldots,l_r>0\\ l_1+\ldots+l_r=a-1}}
 \frac{1}{(1-q_{{\bf n}_{m-1}}^{-n_ms+n_m-a+1+l_{1}})}\frac{\wh v_{X,\F_q,l_1}^{({\bf n}_{m-1})}\ldots\wh v_{X,\F_q,l_r}^{({\bf n}_{m-1})}}{\prod_{j=1}^{r-1}(1-q_{{\bf n}_{m-1}}^{l_j+l_{j+1}})}\\
 \ea
 $$
 where $\dis{\wh\nu_{X,\F_q;N}^{({\bf n}_{m-1})}:=\prod_{k=1}^{N}\wh\zeta_{X,\F_q}^{\,({\bf n}_{m-1})}(k)}$ with $\,\wh\zeta_{X,\F_q}^{\,({\bf n}_{m-1})}(1):=\Res^{~}_{T_{{\bf n}_{m-1}}=1}\wh Z_{X,\F_q}^{\,({\bf n}_{m-1})}(T_{{\bf n}_{m-1}})$.
\end{defn}
In the sequel, for our own convenience, we often write
\be
\wh\zeta_{X,\F_q}^{\,({\bf n}_m)}(s)=\wh\zeta_{X,\F_q}^{\,(n_0,n_1,\ldots,n_m)}(s),\qquad
\wh Z_{X,\F_q}^{\,({\bf n}_m)}(T_{{\bf n}_m}):=\wh\zeta_{X,\F_q}^{\,({\bf n}_m)}(s)\qan {\bf n}_m-j=(n_0,n_1,\ldots,n_m-j),
\ee
and call $\wh Z_{X,\F_q}^{\,({\bf n}_m)}(T_{{\bf n}_m})$ the ${\bf n}_m$-{\it derived Zeta function
of $X$ over $\F_q$}.
\begin{exam}\label{exam1.2} When ${\bf n}_m=({\bf n}_{m-1},1)$ with $n_m=1$, we have, in the definition of 
$\wh\zeta_{X,\F_q}^{\,({\bf n}_m)}(s)$, the most outer summation $\sum_{a=1}^{n_m}=\sum_{a=1}^1$ consisting of a single term $a=1$. This implies that for the second level, the first subsummation
$\sum_{\substack{k_1,\ldots,k_p>0\\ k_1+\ldots+k_p={n_m-a}}}$ is simply
$\sum_{\substack{k_1,\ldots,k_p>0\\ k_1+\ldots+k_p=0}}$, and  hence the corresponding
summand $\equiv 1$ degenerates; similarly, the second
subsummation
$\sum_{\substack{l_1,\ldots,l_r>0\\ l_1+\ldots+l_r=a-1}}$ becomes $\sum_{\substack{l_1,\ldots,l_r>0\\ l_1+\ldots+l_r=0}}$ and hence the corresponding
summand $\equiv 1$ degenerates as well. In addition, for the case under discussion,
$$\wh\zeta_{X,\F_q}^{\,({\bf n}_{m-1})}(n_ms-n_m+a)=\wh\zeta_{X,\F_q}^{\,({\bf n}_{m-1})}(s).
 $$
Therefore, 
 \be
 \wh\zeta_{X,\F_q}^{\,({\bf n}_{m-1},1)}(s)=\wh\zeta_{X,\F_q}^{\,({\bf n}_{m-1})}(s).
 \ee
\end{exam}
We end this subsection with the following comments. If $\wh\zeta_{X,\F_q}^{\,({\bf n}_{m-1})}(s)$ is replaced by $c\wh\zeta_{X,\F_q}^{\,({\bf n}_{m-1})}(s)$ for a certain constant factor $c$. Then
all the special values $\wh\zeta_{X,\F_q}^{\,({\bf n}_{m-1})}(k)$ are replaced by $c\wh\zeta_{X,\F_q}^{\,({\bf n}_{m-1})}(k)$, and all the $\wh v_{X,\F_q,k}^{({\bf n}_{m-1})}$ are replaced by
$c^k\wh v_{X,\F_q,k}^{({\bf n}_{m-1})}$. As a result, using the inductive definition, the ${\bf n}_m$-derived zeta function 
$\wh\zeta_{X,\F_q}^{\,({\bf n}_m)}(s)$ is replaced by $c^{(n_m-a)+1+(a-1)}\wh\zeta_{X,\F_q}^{\,({\bf n}_m)}(s)$, or the same by $c^{n_m}\cdot \wh\zeta_{X,\F_q}^{\,({\bf n}_m)}(s)$.

\section{Standard Zeta Properties of ${\bf n}_m$-Derived Zeta Functions}
In this section, we will establish standard zeta properties for the ${\bf n}_m$-derived zeta function $\wh\zeta_{X,\F_q}^{\,({\bf n}_m)}(s)$ of curves $X$ over $\F_q$. In particular, we show that all these derived zeta functions satisfy the standard functional equation, and are  ratios of  degree 2g polynomials $P_{X,\F_q}^{({\bf n}_m)}(T_{{\bf n}_m})$ by $(1-T_{{\bf n}_m})(1-q_{{\bf n}_m}T_{{\bf n}_m})T_{{\bf n}_m}^{g-1}$.

\subsection{Rationality}

The first zeta property for the ${\bf n}_m$-derived zeta function $\wh\zeta_{X,\F_q}^{\,({\bf n}_m)}(s)$ of curves $X$ over $\F_q$ is the rationality.  Directly, from the inductive definition, we see that for each $(m+1)$-tuple ${\bf n}_m=(n_0,n_1,\ldots, n_m)$ of positive integers, $\wh\zeta_{X,\F_q}^{\,({\bf n}_m)}(s)$ is a rational function of $t:=q^{-s}$. In fact much better can be established.

\begin{thm}[Rationality]\label{thm2.1} Let $X$ be an integral regular projective curve of genus $g$ over $\F_q$.  For each  fixed $(m+1)$-tuple ${\bf n}_m=(n_0,n_1,\ldots, n_m)$ of positive integers, the ${\bf n}_m$-derived Zeta function
$\wh Z_{X,\F_q}^{\,({\bf n}_m)}(T_{{\bf n}_m})$ of $X$ over $\F_q$ is a rational function of $T_{{\bf n}_m}=q^{-n_0n_1\cdots n_m\,s}$.
\end{thm}
\bp
We prove this theorem using an induction on $m$.

When $m=0$, this is the rationality statement of Theorem\,\ref{thm1.1}.

Assume now that the  ${\bf n}_{m-1}$-derived zeta function $\wh Z_{X,\F_q}^{\,({\bf n}_{m-1})}(T_{{\bf n}_{m-1}})$ is a rational function of $T_{{\bf n}_{m-1}}=q^{-n_0n_1\cdots n_{m-1}\,s}$. Then, $\wh Z_{X,\F_q}^{\,({\bf n}_{m-1})}(n_ms-n_m+a)$ is a rational function of 
$$
q^{-n_0n_1\cdots n_{m-1}\,(n_ms)}=T_{{\bf n}_{m-1}}^{n_m}=T_{{\bf n}_m}.
$$ 
This, together with the fact that both
$$
\frac{1}{\Big(1-q_{{\bf n}_{m-1}}^{n_ms-n_m+a+k_{p}}\Big)}\qqan
 \frac{1}{\Big(1-q_{{\bf n}_{m-1}}^{-n_ms+n_m-a+1+l_{1}}\Big)}
 $$
 are rational function of 
 $$
 q_{{\bf n}_{m-1}}^{-n_ms}=\Big(q^{n_0n_1\cdots n_{m-1}}\Big)^{-n_ms}=T_{{\bf n}_{m-1}}^{n_m}=T_{{\bf n}_m}.$$
 Therefore, $\wh Z_{X,\F_q}^{\,({\bf n}_m)}(T_{{\bf n}_m})$ of $X$ over $\F_q$ is a rational function of $T_{{\bf n}_m}$ as well.
\ep

\subsection{Functional equation}\label{sec2.2}
In this subsection, we prove the following important zeta fact for ${\bf n}_{m}$-derived zeta function $\wh \zeta_{X,\F_q}^{\,({\bf n}_{m})}(s)$.
\begin{thm}[Functional Equation]\label{thm2.2}
Let $X$ be an integral regular projective curve of genus $g$ over $\F_q$.  For each  fixed $(m+1)$-tuple ${\bf n}_m=(n_0,n_1,\ldots, n_m)$ of positive integers, the ${\bf n}_m$-derived zeta function
$\wh \zeta_{X,\F_q}^{\,({\bf n}_m)}(s)$ of $X$ over $\F_q$ satisfies thye following standard functional equality
\be
\wh \zeta_{X,\F_q}^{\,({\bf n}_m)}(1-s)=\wh \zeta_{X,\F_q}^{\,({\bf n}_m)}(s)
\ee
\end{thm}
\bp
We prove this theorem using an induction on $m$. When $m=0$, this is the functional equation part of Theorem\,\ref{thm1.1}.

Assume now that the  ${\bf n}_{m-1}$-derived zeta function $\wh \zeta_{X,\F_q}^{\,({\bf n}_{m-1})}(s)$ satisfies the functional function $
\wh \zeta_{X,\F_q}^{\,({\bf n}_{m-1})}(1-s)=\wh \zeta_{X,\F_q}^{\,({\bf n}_{m-1})}(s). $
 Then, for the  ${\bf n}_{m}$-derived zeta function $\wh \zeta_{X,\F_q}^{\,({\bf n}_{m})}(s)$, we have, from Definition\,\ref{defn1.1}(2),
$$\ba 
&\wh \zeta_{X,\F_q}^{\,({\bf n}_{m})}(1-s)\\
=&
q_{{\bf n}_{m-1}}^{\binom{n^{~}_m}{2}(g-1)}\sum_{a=1}^{n_m}\sum_{\substack{k_1,\ldots,k_p>0\\ k_1+\ldots+k_p=n_m-a}}\frac{\wh v_{X,\F_q,k_1}^{({\bf n}_{m-1})}\ldots\wh v_{X,\F_q,k_p}^{({\bf n}_{m-1})}}{\prod_{j=1}^{p-1}(1-q_{{\bf n}_{m-1}}^{k_j+k_{j+1}})} \frac{1}{(1-q_{{\bf n}_{m-1}}^{n_m(1-s)-n_m+a+k_{p}})}\\
 &\hskip 0.20cm\times\wh\zeta_{X,\F_q}^{\,({\bf n}_{m-1})}(n_m(1-s)-n_m+a)\sum_{\substack{l_1,\ldots,l_r>0\\ l_1+\ldots+l_r=a-1}}
 \frac{1}{(1-q_{{\bf n}_{m-1}}^{-n_m(1-s)+n_m-a+1+l_{1}})}\frac{\wh v_{X,\F_q,l_1}^{({\bf n}_{m-1})}\ldots\wh v_{X,\F_q,l_r}^{({\bf n}_{m-1})}}{\prod_{j=1}^{r-1}(1-q_{{\bf n}_{m-1}}^{l_j+l_{j+1}})}\\
=&
q_{{\bf n}_{m-1}}^{\binom{n^{~}_m}{2}(g-1)}\sum_{a=1}^{n_m}\sum_{\substack{k_1,\ldots,k_p>0\\ k_1+\ldots+k_p=n_m-a}}\frac{\wh v_{X,\F_q,k_1}^{({\bf n}_{m-1})}\ldots\wh v_{X,\F_q,k_p}^{({\bf n}_{m-1})}}{\prod_{j=1}^{p-1}(1-q_{{\bf n}_{m-1}}^{k_j+k_{j+1}})} \frac{1}{(1-q_{{\bf n}_{m-1}}^{-n_ms+a+k_{p}})}\\
 &\hskip 0.20cm\times\wh\zeta_{X,\F_q}^{\,({\bf n}_{m-1})}(n_ms-a+1)\sum_{\substack{l_1,\ldots,l_r>0\\ l_1+\ldots+l_r=a-1}}
 \frac{1}{(1-q_{{\bf n}_{m-1}}^{n_ms-a+1+l_{1}})}\frac{\wh v_{X,\F_q,l_1}^{({\bf n}_{m-1})}\ldots\wh v_{X,\F_q,l_r}^{({\bf n}_{m-1})}}{\prod_{j=1}^{r-1}(1-q_{{\bf n}_{m-1}}^{l_j+l_{j+1}})}\\
 &\hskip 1.0cm {\rm(by\ the\ inductive\ assumption\ on\ the\  functional\ equation\ for\ \wh\zeta_{X,\F_q}^{\,({\bf n}_{m-1})}(s))}\\
  =&
q_{{\bf n}_{m-1}}^{\binom{n^{~}_m}{2}(g-1)}\sum_{a=1}^{n_m}\sum_{\substack{k_1,\ldots,k_p>0\\ k_1+\ldots+k_p=a-1}}\frac{\wh v_{X,\F_q,k_1}^{({\bf n}_{m-1})}\ldots\wh v_{X,\F_q,k_p}^{({\bf n}_{m-1})}}{\prod_{j=1}^{p-1}(1-q_{{\bf n}_{m-1}}^{k_j+k_{j+1}})} \frac{1}{(1-q_{{\bf n}_{m-1}}^{-n_ms+n_m-a+1+k_{p}})}\\
 &\hskip 0.20cm\times\wh\zeta_{X,\F_q}^{\,({\bf n}_{m-1})}(n_ms-n_m+a)\sum_{\substack{l_1,\ldots,l_r>0\\ l_1+\ldots+l_r=n_m-a}}
 \frac{1}{(1-q_{{\bf n}_{m-1}}^{n_ms-n_m+a+l_{1}})}\frac{\wh v_{X,\F_q,l_1}^{({\bf n}_{m-1})}\ldots\wh v_{X,\F_q,l_r}^{({\bf n}_{m-1})}}{\prod_{j=1}^{r-1}(1-q_{{\bf n}_{m-1}}^{l_j+l_{j+1}})}\\
 &\hskip 5.0cm {\rm(make\ the\ change\ a\lra n_m-a+1)}\\\ea$$
 $$\ba
=&
q_{{\bf n}_{m-1}}^{\binom{n^{~}_m}{2}(g-1)}\sum_{a=1}^{n_m}\sum_{\substack{k_1,\ldots,k_p>0\\ k_1+\ldots+k_p=n_m-a}}\frac{\wh v_{X,\F_q,k_1}^{({\bf n}_{m-1})}\ldots\wh v_{X,\F_q,k_p}^{({\bf n}_{m-1})}}{\prod_{j=1}^{p-1}(1-q_{{\bf n}_{m-1}}^{k_j+k_{j+1}})} \frac{1}{(1-q_{{\bf n}_{m-1}}^{n_ms-n_m+a+k_{p}})}\\
 &\hskip 0.20cm\times\wh\zeta_{X,\F_q}^{\,({\bf n}_{m-1})}(n_ms-n_m+a)\sum_{\substack{l_1,\ldots,l_r>0\\ l_1+\ldots+l_r=a-1}}
 \frac{1}{(1-q_{{\bf n}_{m-1}}^{-n_ms+n_m-a+1+l_{1}})}\frac{\wh v_{X,\F_q,l_1}^{({\bf n}_{m-1})}\ldots\wh v_{X,\F_q,l_r}^{({\bf n}_{m-1})}}{\prod_{j=1}^{r-1}(1-q_{{\bf n}_{m-1}}^{l_j+l_{j+1}})}\\
 &\hskip 3.0cm {\rm(make\ the\ interchange\ (k_1,\ldots,k_p)\llra (l_1,\ldots,l_r))}\\
=&\wh \zeta_{X,\F_q}^{\,({\bf n}_{m})}(s)\\
\ea$$
as wanted.
\ep
We end this discussion by pointing out that the symmetric exposed in this proof, particularly, in the final  sequences of identities will be used in the next subsection to analyze the singularities of the derived zeta functions.

\subsection{Singularities}\label{sec2.3}
Before we state the main result of this section, let us examine the structure of $\wh \zeta_{X,\F_q}^{\,({\bf n}_{m})}(s)$ in more details.
Set, for $1\leq a\leq n_m$,
$$\ba
&\wh \zeta_{X,\F_q}^{\,({\bf n}_{m}),[a]}(s):=\wh Z_{X,\F_q}^{\,({\bf n}_{m}),[a]}(T_{{\bf n}_m})\\
:=&q_{{\bf n}_{m-1}}^{\binom{n^{~}_m}{2}(g-1)}\sum_{\substack{k_1,\ldots,k_p>0\\ k_1+\ldots+k_p=n_m-a}}\frac{\wh v_{X,\F_q,k_1}^{({\bf n}_{m-1})}\ldots\wh v_{X,\F_q,k_p}^{({\bf n}_{m-1})}}{\prod_{j=1}^{p-1}(1-q_{{\bf n}_{m-1}}^{k_j+k_{j+1}})} \frac{1}{(1-q_{{\bf n}_{m-1}}^{n_ms-n_m+a+k_{p}})}\\
 &\hskip 0.20cm\times\wh\zeta_{X,\F_q}^{\,({\bf n}_{m-1})}(n_ms-n_m+a)\sum_{\substack{l_1,\ldots,l_r>0\\ l_1+\ldots+l_r=a-1}}
 \frac{1}{(1-q_{{\bf n}_{m-1}}^{-n_ms+n_m-a+1+l_{1}})}\frac{\wh v_{X,\F_q,l_1}^{({\bf n}_{m-1})}\ldots\wh v_{X,\F_q,l_r}^{({\bf n}_{m-1})}}{\prod_{j=1}^{r-1}(1-q_{{\bf n}_{m-1}}^{l_j+l_{j+1}})}\\
\ea$$
Then \be
\wh \zeta_{X,\F_q}^{\,({\bf n}_{m})}(s)=\sum_{a=1}^{n_m}\wh \zeta_{X,\F_q}^{\,({\bf n}_{m}),[a]}(s)
\ee
\begin{defn}\label{defn2.1} Let $X$ be an integral regular projective curve of genus $g$ over $\F_q$.  For each  fixed $(m+1)$-tuple ${\bf n}_m=(n_0,n_1,\ldots, n_m)$ of positive integers, we define
\begin{enumerate}
\item [(1)]  For any $1\leq n\leq n_m$, the  rational function
$\De_{X,\F_q;{\bf n}_{m};n}(T_{{\bf n}_m})$,  resp. the polynomial $\Ga_{X,\F_q;{\bf n}_{m};n}(T_{{\bf n}_m})$,  of $T_{{\bf n}_m}$ is defined by
$$\ba
\De_{X,\F_q;{\bf n}_{m};n}(T_{{\bf n}_m}):
=&
\sum_{\substack{k_1,\ldots,k_p>0\\ k_1+\ldots+k_p=n}}\frac{\wh v_{X,\F_q,k_1}^{({\bf n}_{m-1})}\ldots\wh v_{X,\F_q,k_p}^{({\bf n}_{m-1})}}{\prod_{j=1}^{p-1}\Big(1-q_{{\bf n}_{m-1}}^{k_j+k_{j+1}}\Big)} \frac{1}{\Big(q_{{\bf n}_{m-1}}^{k_{p}}T_{{\bf n}_{m}}-1\Big)}
\ea
$$
resp.
\be
\Ga_{X,\F_q;{\bf n}_{m},n}(T_{{\bf n}_m}):=\De_{X,\F_q;{\bf n}_{m-1},n}(T_{{\bf n}_m})\cdot\prod_{\substack{\ell=1}}^n\Big(q_{{\bf n}_{m-1}}^{\ell} T_{{\bf n}_{m}}-1\Big)
\ee
\item[(2)] For any positive integer $n\geq 1$, the {\it $n$-th ${\bf n}_m$-derived $\b$-invariant for $X$ over $\F_q$} by
\be\ba
\b_{X,\F_q;{\bf n}_{m-1};n}:=&q_{{\bf n}_{m-1}}^{\binom{n}{2}(g-1)}\sum_{\substack{k_1,\ldots,k_p>0\\ k_1+\ldots+k_p=n}}\frac{\wh v_{X,\F_q,k_1}^{({\bf n}_{m-1})}\ldots\wh v_{X,\F_q,k_p}^{({\bf n}_{m-1})}}{\prod_{j=1}^{p-1}\Big(1-q_{{\bf n}_{m-1}}^{k_j+k_{j+1}}\Big)}\\
 \b_{X,\F_q;{\bf n}_{m-1};n}'':=&q_{{\bf n}_{m-1}}^{-\binom{n}{2}(g-1)}\b_{X,\F_q;{\bf n}_{m-1};n}
 =\sum_{\substack{k_1,\ldots,k_p>0\\ k_1+\ldots+k_p=n}}\frac{\wh v_{X,\F_q,k_1}^{({\bf n}_{m-1})}\ldots\wh v_{X,\F_q,k_p}^{({\bf n}_{m-1})}}{\prod_{j=1}^{p-1}\Big(1-q_{{\bf n}_{m-1}}^{k_j+k_{j+1}}\Big)}.
\ea\ee
\end{enumerate}
\end{defn}
For our own convenience, we also set
\be
\b_{X,\F_q}^{({\bf n}_m)}:=\b_{X,\F_q;{\bf n}_m}:=\b_{X,\F_q;{\bf n}_{m-1};n_m}
\ee
and denote the leading coefficient and the constant term of the polynomial $\Ga_{X,\F_q;{\bf n}_{m};n}(T_{{\bf n}_m})$ in by
$L_{X,\F_q;{\bf n}_{m};n}$ and $C_{X,\F_q;{\bf n}_{m};n}$, respectively.

\begin{lem}\label{lem2.3} With the same notation as above, we have
\begin{enumerate}
\item [(1)] $L_{X,\F_q;{\bf n}_{m};n}$ and $C_{X,\F_q;{\bf n}_{m};n}$ are given by
$$\ba
q_{{\bf n}_{m-1}}^{-\binom{n+1}{2}}L_{X,\F_q;{\bf n}_{m};n}=&\sum_{\substack{k_1,\ldots,k_p>0\\ k_1+\ldots+k_p=n}}\frac{\wh v_{X,\F_q,k_1}^{({\bf n}_{m-1})}\ldots\wh v_{X,\F_q,k_p}^{({\bf n}_{m-1})}}{\prod_{j=1}^{p-1}\Big(1-q_{{\bf n}_{m-1}}^{k_j+k_{j+1}}\Big)}q_{{\bf n}_{m-1}}^{-k_{p}}\\
(-1)^nC_{X,\F_q;{\bf n}_{m};n}=&\b_{X,\F_q;{\bf n}_{m-1};n}''=\sum_{\substack{k_1,\ldots,k_p>0\\ k_1+\ldots+k_p=n}}\frac{\wh v_{X,\F_q,k_1}^{({\bf n}_{m-1})}\ldots\wh v_{X,\F_q,k_p}^{({\bf n}_{m-1})}}{\prod_{j=1}^{p-1}\Big(1-q_{{\bf n}_{m-1}}^{k_j+k_{j+1}}\Big)}\\
\ea$$
\item[(2)]  $\Ga_{X,\F_q;{\bf n}_{m};n}(T_{{\bf n}_m})$ is a polynomial in $T_{{\bf n}_m}$ of degree $(n-1)$, provided that $L_{X,\F_q;{\bf n}_{m};n}$ is non-zero. In addition,
 $\dis{\frac{\De_{X,\F_q;{\bf n}_{m-1},n}(q_{{\bf n}_{m-1}}^{-n}T_{{\bf n}_m}^{-1})\cdot\prod_{\substack{\ell=0}}^n\Big(q_{{\bf n}_{m-1}}^{\ell} T_{{\bf n}_{m}}-1\Big)}{\ T_{{\bf n}_m}}}$ is a degree $n$ polynomial of $T_{{\bf n}_m}$ with leading coefficient
$-q_{{\bf n}_{m-1}}^{\binom{n+1}{2}}\b_{X,\F_q;{\bf n}_{m-1};n}''.$
\item [(3)] 
$\dis{
\wh Z_{X,\F_q}^{\,({\bf n}_{m}),[a]}(T_{{\bf n}_m})=\De_{X,\F_q;{\bf n}_{m};n_m-a}(q_{{\bf n}_{m-1}}^{-n_m+a}T_{{\bf n}_m}^{-1})\cdot 
\wh\zeta_{X,\F_q}^{\,({\bf n}_{m})}(n_ms-n_m+a)\cdot\De_{X,\F_q;{\bf n}_{m};a-1}(q_{{\bf n}_{m-1}}^{n_m-a+1}T_{{\bf n}_m})}$
\end{enumerate}
\end{lem}
\bp
All the statements comes directly from the definition. As an illustration, we give a proof of the second assertion of (2). By Definition\,\ref{defn2.1}(1),
$$\ba
&\frac{1}{\ T_{{\bf n}_m}}\De_{X,\F_q;{\bf n}_{m};n}(q_{{\bf n}_{m-1}}^{-n}T_{{\bf n}_m}^{-1})\cdot \prod_{\ell=0}^n\Big(q_{{\bf n}_{m-1}}^{\ell}T_{{\bf n}_m}-1\Big)\\
=&\frac{1}{\ T_{{\bf n}_m}}\cdot \prod_{\ell=0}^n\Big(q_{{\bf n}_{m-1}}^{n-\ell}T_{{\bf n}_m}-1\Big) 
\sum_{\substack{k_1,\ldots,k_p>0\\ k_1+\ldots+k_p=n}}\frac{\wh v_{X,\F_q,k_1}^{({\bf n}_{m-1})}\ldots\wh v_{X,\F_q,k_p}^{({\bf n}_{m-1})}}{\prod_{j=1}^{p-1}\Big(1-q_{{\bf n}_{m-1}}^{k_j+k_{j+1}}\Big)} \frac{1}{\Big(q_{{\bf n}_{m-1}}^{-n+k_{p}}T_{{\bf n}_{m}}^{-1}-1\Big)}\\
=&\frac{1}{\ T_{{\bf n}_m}}\cdot \prod_{\ell=0}^n\Big(q_{{\bf n}_{m-1}}^{n-\ell}T_{{\bf n}_m}-1\Big) 
\sum_{\substack{k_1,\ldots,k_p>0\\ k_1+\ldots+k_p=n}}\frac{\wh v_{X,\F_q,k_1}^{({\bf n}_{m-1})}\ldots\wh v_{X,\F_q,k_p}^{({\bf n}_{m-1})}}{\prod_{j=1}^{p-1}\Big(1-q_{{\bf n}_{m-1}}^{k_j+k_{j+1}}\Big)} \frac{q_{{\bf n}_{m-1}}^{n-k_{p}}T_{{\bf n}_{m}}}{\Big(1-q_{{\bf n}_{m-1}}^{n-k_{p}}T_{{\bf n}_{m}}\Big)}
\ea$$
is a degree $n$ polynomial in $T_{{\bf n}_m}$, whose leading coefficient is given by
$$
-\sum_{\substack{k_1,\ldots,k_p>0\\ k_1+\ldots+k_p=n}}\frac{\wh v_{X,\F_q,k_1}^{({\bf n}_{m-1})}\ldots\wh v_{X,\F_q,k_p}^{({\bf n}_{m-1})}}{\prod_{j=1}^{p-1}\Big(1-q_{{\bf n}_{m-1}}^{k_j+k_{j+1}}\Big)}q_{{\bf n}_{m-1}}^{\sum_{\ell=0}^n(n-\ell)}=-q_{{\bf n}_{m-1}}^{\binom{n+1}{2}}\sum_{\substack{k_1,\ldots,k_p>0\\ k_1+\ldots+k_p=n}}\frac{\wh v_{X,\F_q,k_1}^{({\bf n}_{m-1})}\ldots\wh v_{X,\F_q,k_p}^{({\bf n}_{m-1})}}{\prod_{j=1}^{p-1}\Big(1-q_{{\bf n}_{m-1}}^{k_j+k_{j+1}}\Big)},
$$
as wanted.
\ep
\begin{rmk}
From this lemma, for a fixed $1\leq a\leq n_m$,  besides the singularities coming from
the ${\bf n}_{m-1}$-derived Zeta function $\wh Z_{X,\F_q}^{\,({\bf n}_{m})}(T_{{\bf n}_{m-1}})$, we conclude that 
$\wh Z_{X,\F_q}^{\,({\bf n}_{m}),[a]}(T_{{\bf n}_m})$ admits at least singularities at $T_{{\bf n}_m}=q_{{\bf n}_{m-1}}^{-n_m+a+k_{p}}$
for $k_p=1,\ldots n_m-a$ and at $q_{{\bf n}_{m-1}}^{-n_m+a-1-l_{1}}$ for $l_1=1,\ldots, a-1$, or equivalently at
$$
T_{{\bf n}_m}=q_{{\bf n}_{m-1}}^{-n_m},q_{{\bf n}_{m-1}}^{-n_m+1},\ldots, q_{{\bf n}_{m-1}}^{-n_m+a-2},q_{{\bf n}_{m-1}}^{-n_m+a+1}, q_{{\bf n}_{m-1}}^{-n_m+a+2},\ldots, q_{{\bf n}_{m-1}}^{0}=1.
$$
\end{rmk}
Accordingly, it appears that the ${\bf n}_{m}$-derived Zeta function $\wh Z_{X,\F_q}^{\,({\bf n}_{m})}(T_{{\bf n}_m})$ might admit many singularities and could be hardly handled. 
However,  one beauty of the derived zeta function  is that there are perfect cancelations among these singularities of the $\wh Z_{X,\F_q}^{\,({\bf n}_{m}),[a]}(T_{{\bf n}_m})$'s. 

\begin{thm}[Singularities] \label{thm2.4} Let $X$ be an integral regular projective curve of genus $g$ over $\F_q$.  For each  fixed $(m+1)$-tuple ${\bf n}_m=(n_0,n_1,\ldots, n_m)$ of positive integers, 
\begin{enumerate}
\item [(1)] The ${\bf n}_m$-derived Zeta function $Z_{X,\F_q}^{\,({\bf n}_{m})}(T_{{\bf n}_m}):=T_{{\bf n}_m}^{g-1}\cdot\wh Z_{X,\F_q}^{\,({\bf n}_{m})}(T_{{\bf n}_m})$ admits only two singularities, namely two simple poles at $T=1$ and $T=1/q_{{\bf n}_m}$ whose residues are given by
\be
\Res_{T=1}\wh Z_{X,\F_q}^{\,({\bf n}_{m})}(T_{{\bf n}_m})=-\Res_{T=1/q_{{\bf n}_m}}\wh Z_{X,\F_q}^{\,({\bf n}_{m})}(T_{{\bf n}_m})=\b_{X,\F_q;{\bf n}_{m-1};n_m}=\b_{X,\F_q;{\bf n}_m}.
\ee
Furthermore,
\item[(2)] As a  rational function of $T_{{\bf n}_m}$, there exists a degree $2g$ polynomial $P_{X,\F_q}^{\,({\bf n}_{m})}(T_{{\bf n}_m})$ with rational coefficients such that
\be
\wh Z_{X,\F_q}^{\,({\bf n}_{m})}(T_{{\bf n}_m})=\frac{P_{X,\F_q}^{\,({\bf n}_{m})}(T_{{\bf n}_m})}{\big(1-T_{{\bf n}_m}\big)\big(1-q_{{\bf n}_m}T_{{\bf n}_m}\big)T_{{\bf n}_m}^{g-1}}
\ee
\end{enumerate}
\end{thm}
\bp
We prove this theorem using an induction on $m$.

When $m=0$, (1), resp. (2) is the Singularities, resp. Rationality, part of Theorem\,\ref{thm1.1} for rank $n$-zeta non-abelian functions, by using the special uniformity of zetas, i.e. Theorem\,\ref{thm1.3}.

Assume now that the assertions (1) and (2) in the theorem hold for ${\bf n}_{m-1}$-derived Zeta function $\wh Z_{X,\F_q}^{\,({\bf n}_{m-1})}(T_{{\bf n}_{m-1}})$. We prove (1) and (2) hold for ${\bf n}_{m}$-derived Zeta function $\wh Z_{X,\F_q}^{\,({\bf n}_{m})}(T_{{\bf n}_m})$.

By the remark immediately after Lemma\,\ref{lem2.3} and the inductive hypothesis, the possible singularities of $\wh Z_{X,\F_q}^{\,({\bf n}_{m-1})}(T_{{\bf n}_{m-1}})$ are all simple poles located at $q_{{\bf n}_{m-1}}^{-\ell}$ where $0\leq \ell\leq n_m$. Our first task is to show that except for $T_{{\bf n}_{m}}=q_{{\bf n}_{m-1}}^{-0}=1$ and $q_{{\bf n}_{m}}^{-1}$, the derived zeta $\wh Z_{X,\F_q}^{\,({\bf n}_{m-1})}(T_{{\bf n}_{m-1}})$ is holomorphic at all $T_{{\bf n}_{m}}=q_{{\bf n}_{m-1}}^{-\ell}$ for $1\leq \ell\leq n_m-1$. Indeed, this can be verified as what we have done in the final section of \cite{WZ2}. To see this, set
$$
R_{\ell,a}:=\Res_{T_{{\bf n}_{m}}=q_{{\bf n}_{m-1}}^{-\ell}}q_{{\bf n}_{m-1}}^{-\binom{n^{~}_m}{2}(g-1)}\wh Z_{X,\F_q}^{\,({\bf n}_{m}),[a]}(T_{{\bf n}_{m}})\qan R_{\ell}:=\Res_{T_{{\bf n}_{m}}=q_{{\bf n}_{m-1}}^{-\ell}}q_{{\bf n}_{m-1}}^{-\binom{n^{~}_m}{2}(g-1)}\wh Z_{X,\F_q}^{\,({\bf n}_{m})}(T_{{\bf n}_{m}}).
$$
Then
$R_{\ell}=\sum_{a=1}^{n_m}R_{\ell,a}$.
Furthermore,
recall that there is a natural  symmetry $(a,T_{{\bf n}_{m}})\llra(n_m-a-1, q_{{\bf n}_{m}}T_{{\bf n}_{m}}^{-1})$ for the summands of $\wh Z_{X,\F_q}^{\,({\bf n}_{m})}(T_{{\bf n}_{m}})$
as exposed in our proof of the functional equation of Theorem\,\ref{thm3.2} (see also the comments at the end of \S\ref{sec2.2}). Hence
$$
R_{\ell,a}=-R_{n_m-\ell,n_m-a-1}\qqan R_\ell=S_\ell-S_{n_m-\ell}
$$
where
$$
S_\ell=\sum_{a=0}^{\ell-1}R_{\ell,a}
$$
In addition, directly from Lemma\,\ref{lem2.3}(3), for $0\leq a\leq \ell-1$, we have
$$\ba
 R_{\ell,a}=&\De_{X,\F_q;{\bf n}_{m};n_m-a}(q_{{\bf n}_{m-1}}^{\ell-a})\cdot \wh Z_{X,\F_q}^{\,({\bf n}_{m-1})}(q_{{\bf n}_{m-1}}^{a-\ell})\wh v_{X,\F_q,\ell-a-1}^{\,({\bf n}_{m-1})}\cdot \De_{X,\F_q;{\bf n}_{m};a-1}(q_{{\bf n}_{m-1}}^{\ell-a-1})\\
 =&\De_{X,\F_q;{\bf n}_{m};n_m-a}(q_{{\bf n}_{m-1}}^{\ell-a})\cdot \wh v_{X,\F_q,\ell-a}^{\,({\bf n}_{m-1})}\cdot \De_{X,\F_q;{\bf n}_{m};a-1}(q_{{\bf n}_{m-1}}^{\ell-a-1})\\
 &\hskip 4.0cm\Big({\rm since}\ \wh v_{X,\F_q,\ell-a}^{\,({\bf n}_{m-1})}=\wh Z_{X,\F_q}^{\,({\bf n}_{m-1})}(q_{{\bf n}_{m-1}}^{a-\ell})\cdot\wh v_{X,\F_q,\ell-a-1}^{\,({\bf n}_{m-1})}\Big)\\
 =&\sum_{1\leq s<p\leq n_m}\sum_{\substack{k_1,\ldots,k_p\geq 1\\ k_1+\ldots+k_p=n_m\\ k_1+\ldots+k_{s-1}=a, n_s=\ell-a}}\frac{\wh v_{X,\F_q,k_1}^{({\bf n}_{m-1})}\ldots\wh v_{X,\F_q,k_p}^{({\bf n}_{m-1})}}{\prod_{j=1}^{p-1}\Big(1-q_{{\bf n}_{m-1}}^{k_j+k_{j+1}}\Big)}.
\ea
$$
Consequently,
$$
S_{\ell}=\sum_{1\leq s<p\leq n_m}\sum_{\substack{k_1,\ldots,k_p\geq 1\\ k_1+\ldots+k_p=n_m\\ k_1+\ldots+k_{s}=\ell}}\frac{\wh v_{X,\F_q,k_1}^{({\bf n}_{m-1})}\ldots\wh v_{X,\F_q,k_p}^{({\bf n}_{m-1})}}{\prod_{j=1}^{p-1}\Big(1-q_{{\bf n}_{m-1}}^{k_j+k_{j+1}}\Big)}.
$$
which is obviously symmetric under $\ell\llra n_m-\ell$ by replacing $k_j$ by $k_{p+1-j}$. Therefore,
$$
R_\ell=S_\ell-S_{n_m-\ell}=0.
$$
We claim that with a similar argument, as to be seen in the next subsection, we have that
$$
R_n=-R_0=\sum_{1\leq s<p\leq n_m}\sum_{\substack{k_1,\ldots,k_p\geq 1\\ k_1+\ldots+k_p=n_m\\ }}\frac{\wh v_{X,\F_q,k_1}^{({\bf n}_{m-1})}\ldots\wh v_{X,\F_q,k_p}^{({\bf n}_{m-1})}}{\prod_{j=1}^{p-1}\Big(1-q_{{\bf n}_{m-1}}^{k_j+k_{j+1}}\Big)}=\b_{X,\F_q;{\bf n}_{m-1};n_m}''.
$$
This proves (1) and (2) since the multiple pole at $T_{{\bf n}_{m}}=0$ with multiplicity $g-1$ comes directly and solely from the level $(m-1)$-derived zeta function, by our inductive assumption,
$$
\wh Z_{X,\F_q}^{\,({\bf n}_{m-1})}(T_{{\bf n}_{m-1}})=\frac{P_{X,\F_q}^{\,({\bf n}_{m-1})}(T_{{\bf n}_{m-1}})}{\big(1-T_{{\bf n}_{m-1}}\big)\big(1-q_{{\bf n}_{m-1}}T_{{\bf n}_{m-1}}\big)T_{{\bf n}_{m-1}}^{g-1}}
$$
and the relation $T_{{\bf n}_{m}}=T_{{\bf n}_{m-1}}^{n_m}$. 
\ep
As a direct consequence of Theorems\,\ref{thm2.1},\,\ref{thm2.2} and \ref{thm2.4}, we have the following
\begin{cor}\label{cor2.5} Let $X$ be an integral regular projective curve of genus $g$ on $\F_q$. Then
\be\label{eq39}\ba
\wh \zeta_{X,\F_q}^{({\bf n}_{m})}(s)=\wh Z_{X,\F_q}^{({\bf n}_{m})}(T_{{\bf n}_{m}})
=&
\Big(\sum_{\ell=0}^{g-2}\a_{X,\F_q}^{({\bf n}_{m})}(\ell)\Big(T_{{\bf n}_{m}}^{\ell-(g-1)}+q_{{\bf n}_{m}}^{(g-1)-\ell}T_{{\bf n}_{m}}^{(g-1)-\ell}\Big)
+\a_{X,\F_q}^{({\bf n}_{m})}(g-1))\Big)\\
&\hskip 4.0cm+\frac{(q_{{\bf n}_{m}}-1)T_{{\bf n}_{m}}\b_{X,\F_q;{\bf n}_{m}}}{(1-T_{{\bf n}_{m}})(1-q_{{\bf n}_{m}}T_{{\bf n}_{m}})}\\
\ea\ee
for some ${\bf n}_m$-derived alpha invariants $\Big\{\a_{X,\F_q}^{({\bf n}_{m})}(\ell)\Big\}_{\ell=0}^{g-1}$ and the beta invariant $\b_{X,\F_q;{\bf n}_{m}}$.
of $X/\F_q$. 
\end{cor}
\bp
By Theorem\,\ref{thm2.4}, $\wh \zeta_{X,\F_q}^{({\bf n}_{m})}(s)$ admits only two simple poles at 
$s=0,1$ with residue $\b_{X,\F_q;{\bf n}_{m}}$, we conclude that $\wh Z_{X,\F_q}^{({\bf n}_{m})}(T_{{\bf n}_{m}})$ admits two simples at
$T_{{\bf n}_{m}}=1$ and $q_{{\bf n}_{m}}^{-1}$. This gives the final term of \eqref{eq39}. In addition, by the rationality and the singularities explained above for  $\wh \zeta_{X,\F_q}^{({\bf n}_{m})}(s)$, we conclude that
$\wh Z_{X,\F_q}^{({\bf n}_{m})}(T_{{\bf n}_{m}})$  should take the form
$\frac{P_{X,\F_q}^{({\bf n}_m)}(T_{{\bf n}_{m}})}{(1-T_{{\bf n}_{m}})(1-q_{{\bf n}_{m}}T_{{\bf n}_{m}})T_{{\bf n}_{m}}^a}$
 for a certain polynomial $P_{X,\F_q}^{({\bf n}_m)}(T_{{\bf n}_{m}})$ and a certain non-negative integer $a$, by noticing the fact that in terms of $s$, $\frac{1}{q_{{\bf n}_{m}}^{-as}}$ admits no singularity.
Therefore, finally, by applying the functional equation
\be
\wh \zeta_{X,\F_q}^{({\bf n}_{m})}(1-s)=\wh \zeta_{X,\F_q}^{({\bf n}_{m})}(s)\qquad{\rm
or\ better}\qquad
\wh Z_{X,\F_q}^{({\bf n}_{m})}\Big(\frac{1}{q_{{\bf n}_{m}}T_{{\bf n}_{m}}}\Big)=\wh Z_{X,\F_q}^{({\bf n}_{m})}(T_{{\bf n}_{m}}), 
\ee
we easily conclude that $a=g-1$ and the coefficients of $P_{X,\F_q}^{({\bf n}_m)}(T_{{\bf n}_{m}})$ should take the form in \eqref{eq39}.

With \eqref{eq39} established, what is left is to verify that $\b_{X,\F_q;{\bf n}_{m}}$ is given by the closed formula in Definition\,\ref{defn2.1}(2). This will be done in the following subsection, in particular, in Theorem\,\ref{thm2.6} below.
\ep

\subsection{${\bf n}_m$-Derived beta invariants: Justification of the closed formula}

By definition and the special uniformity of zetas,
$$\ba
&\b_{X,\F_q;{\bf n}_{m}}=\Res_{T_{{\bf n}_{m}}=1}\wh Z_{X,\F_q}^{({\bf n}_{m})}(T_{{\bf n}_{m}})\\
=&q_{{\bf n}_{m-1}}^{\binom{n_m}{2}(g-1)}\Res_{T_{{\bf n}_{m}}=1}\Bigg(\sum_{a=1}^{n_m}\sum_{\substack{k_1,\ldots,k_p>0\\ k_1+\ldots+k_p=n_m-a}}\frac{\wh v^{({\bf n}_{m-1})}_{k_1}\ldots\wh v^{({\bf n}_{m-1})}_{k_p}}{\prod_{j=1}^{p-1}(1-q_{{\bf n}_{m-1}}^{k_j+k_{j+1}})} \frac{1}{(1-q_{{\bf n}_{m-1}}^{n_ms-n_m+a+k_{p}})}\\
 &\hskip 2.0cm\times\wh\zeta_{X,\F_q}^{\,({\bf n}_{m-1})}(n_ms-n_m+a)\sum_{\substack{l_1,\ldots,l_r>0\\ l_1+\ldots+l_r=a-1}}
 \frac{1}{(1-q_{{\bf n}_{m-1}}^{-n_ms+n_m-a+1+l_{1}})}\frac{\wh v^{({\bf n}_{m-1})}_{l_1}\ldots\wh v^{({\bf n}_{m-1})}_{l_r}}{\prod_{j=1}^{r-1}(1-q_{{\bf n}_{m-1}}^{l_j+l_{j+1}})}
\Bigg)\\
=&q_{{\bf n}_{m-1}}^{\binom{n_m}{2}(g-1)}\Res_{T_{{\bf n}_{m}}=1}\Bigg(\sum_{a=1}^{n_m}\sum_{\substack{k_1,\ldots,k_p, l_1,\ldots,l_r>0\\ l_1+\ldots+l_r+ k_1+\ldots+k_p=n_m-1\\
 l_1+\ldots+l_r=a-1,\  k_1+\ldots+k_p=n_m-a
}}\frac{\wh v^{({\bf n}_{m-1})}_{k_1}\ldots\wh v^{({\bf n}_{m-1})}_{k_p}\wh v^{({\bf n}_{m-1})}_{l_1}\ldots\wh v^{({\bf n}_{m-1})}_{l_r}}{\prod_{j=1}^{p-1}(1-q_{{\bf n}_{m-1}}^{k_j+k_{j+1})}\prod_{j=1}^{r-1}(1-q_{{\bf n}_{m-1}}^{l_j+l_{j+1}})} \\
&\hskip 3.0cm\times\frac{1}{(1-q_{{\bf n}_{m-1}}^{-n_ms+n_m-a+1+l_{1}})(1-q_{{\bf n}_{m-1}}^{n_ms-n_m+a+k_{p}})}\wh\zeta_{X,\F_q}^{\,({\bf n}_{m-1})}(n_ms-n_m+a) 
\Bigg)\\
=&q_{{\bf n}_{m-1}}^{\binom{n_m}{2}(g-1)}\Res_{T_{{\bf n}_{m}}=1}\Bigg(\sum_{a=1}^{n_m}\sum_{\substack{k_1,\ldots,k_p, l_1,\ldots,l_r>0\\ l_1+\ldots+l_r+ l_{r+1}+\ldots+k_{r+p}=n_m-1\\
 l_1+\ldots+l_r=a-1,\  l_{r+1}+\ldots+l_{r+p}=n_m-a
}}\frac{\wh v^{({\bf n}_{m-1})}_{l_1}\ldots\wh v^{({\bf n}_{m-1})}_{l_r}\wh v^{({\bf n}_{m-1})}_{l_{r+1}}\ldots\wh v^{({\bf n}_{m-1})}_{l_{r+p}}(1-q_{{\bf n}_{m-1}}^{l_1+l_{r+p}})}{\prod_{j=1}^{r+p-1}(1-q_{{\bf n}_{m-1}}^{l_j+l_{j+1}})} \\
&\hskip 3.0cm\times\frac{1}{(1-q_{{\bf n}_{m-1}}^{-n_ms+n_m-a+1+l_{1}})(1-q_{{\bf n}_{m-1}}^{n_ms-n_m+a+l_{r+p}})}\wh\zeta_{X,\F_q}^{\,({\bf n}_{m-1})}(n_ms-n_m+a) 
\Bigg)\\\ea$$
$$
\ba
=&q_{{\bf n}_{m-1}}^{\binom{n_m}{2}(g-1)}\Res_{T_{{\bf n}_{m}}=1}\Bigg(\sum_{\substack{d_1,\ldots,d_h>0\\ d_1+\ldots+d_h=n_m-1
}}\frac{\wh v^{({\bf n}_{m-1})}_{d_1}\ldots\wh v^{({\bf n}_{m-1})}_{d_h}}{\prod_{j=1}^{h-1}(1-q_{{\bf n}_{m-1}}^{d_j+d_{j+1}})} \\
\times&
\left(\frac{\wh\zeta_{X,\F_q}^{({\bf n}_{m-1})}(n_ms)}{(1-q_{{\bf n}_{m-1}}^{-n_ms+d_1+1})}+
\sum_{i=1}^{h-1}\frac{(1-q_{{\bf n}_{m-1}}^{d_i+d_{i+1}})\cdot\wh\zeta_{X,\F_q}^{({\bf n}_{m-1})}(n_ms-(d_1+\ldots+d_i))}
{(1-q_{{\bf n}_{m-1}}^{n_ms-(d_1+\ldots+d_{i-1})})(1-q_{{\bf n}_{m-1}}^{-n_ms+d_1+\ldots+d_{i+1}+1})}+\frac{\wh\zeta_{X,\F_q}^{({\bf n}_{m-1})}(n_ms-n_m+1)}
{(1-q_{{\bf n}_{m-1}}^{n_ms-(d_1+\ldots+d_{h-1})})} \right)
\Bigg)\\
=&q_{{\bf n}_{m-1}}^{\binom{n_m}{2}(g-1)}\Res_{T_{{\bf n}_{m}}=1}\Bigg(\sum_{\substack{d_1,\ldots,d_h>0\\ d_1+\ldots+d_h=n_m-1
}}\frac{\wh v^{({\bf n}_{m-1})}_{d_1}\ldots\wh v^{({\bf n}_{m-1})}_{d_h}}{\prod_{j=1}^{h-1}(1-q_{{\bf n}_{m-1}}^{d_j+d_{j+1}})} \\&\times
\left(\frac{\wh\zeta_{X,\F_q}^{({\bf n}_{m-1})}(n_ms)}{(1-q_{{\bf n}_{m-1}}^{d_1+1}T_{{\bf n}_{m}})}+
\sum_{i=1}^{h-1}\frac{(1-q_{{\bf n}_{m-1}}^{d_i+d_{i+1}})\cdot q_{{\bf n}_{m-1}}^{d_1+\ldots+d_{i-1}}T_{{\bf n}_{m}}\wh\zeta_{X,\F_q}^{({\bf n}_{m-1})}(n_ms-(d_1+\ldots+d_i))}
{(q_{{\bf n}_{m-1}}^{d_1+\ldots+d_{i-1}}T_{{\bf n}_{m}}-1)(1-q_{{\bf n}_{m-1}}^{d_1+\ldots+d_{i+1}+1}T_{{\bf n}_{m}})}\right.\\
&\hskip 6.0cm\left.+\frac{q_{{\bf n}_{m-1}}^{d_1+\ldots+d_{h-1}}T_{{\bf n}_{m}}\wh\zeta_{X,\F_q}^{({\bf n}_{m-1})}(n_ms-n_m+1)}
{(q_{{\bf n}_{m-1}}^{d_1+\ldots+d_{h-1}}T_{{\bf n}_{m}}-1)} \right)
\Bigg)\\
=&q_{{\bf n}_{m-1}}^{\binom{n_m}{2}(g-1)}\Bigg(\sum_{\substack{d_1,\ldots,d_h>0\\ d_1+\ldots+d_h=n_m-1
}}\frac{\wh v^{({\bf n}_{m-1})}_{d_1}\ldots\wh v^{({\bf n}_{m-1})}_{d_h}}{\prod_{j=1}^{h-1}(1-q_{{\bf n}_{m-1}}^{d_j+d_{j+1}})}\cdot
\Res_{T_{{\bf n}_{m}}=1}\left(\frac{\wh\zeta_{X,\F_q}^{({\bf n}_{m-1})}(n_ms)}{(1-q_{{\bf n}_{m-1}}^{d_1+1}T_{{\bf n}_{m}})}\right.\\
&\hskip 2.0cm+
\sum_{i=1}^{h-1}\frac{(1-q_{{\bf n}_{m-1}}^{d_i+d_{i+1}})\cdot q_{{\bf n}_{m-1}}^{d_1+\ldots+d_{i-1}}T_{{\bf n}_{m}}\wh\zeta_{X,\F_q}^{({\bf n}_{m-1})}(n_ms-(d_1+\ldots+d_i))}
{(q_{{\bf n}_{m-1}}^{d_1+\ldots+d_{i-1}}T_{{\bf n}_{m}}-1)(1-q_{{\bf n}_{m-1}}^{d_1+\ldots+d_{i+1}+1}T_{{\bf n}_{m}})}\\
&\hskip 5.0cm\left.+\frac{q_{{\bf n}_{m-1}}^{d_1+\ldots+d_{h-1}}T_{{\bf n}_{m}}\wh\zeta_{X,\F_q}^{({\bf n}_{m-1})}(n_ms-n_m+1)}
{(q_{{\bf n}_{m-1}}^{d_1+\ldots+d_{h-1}}T_{{\bf n}_{m}}-1)} \right)
\Bigg)\\
\ea$$
Recall that, by Corollary\,\ref{cor2.5}, we have
$$\ba
&\wh\zeta_{X,\F_q}^{({\bf n}_{m-1})}(n_ms)\\
=&\frac{
\Big(\sum_{\ell=0}^{g-2}\a_{X,\F_q}^{({\bf n}_{m-1})}(\ell)\Big(T_{{\bf n}_{m}}^{m}+q_{{\bf n}_{m-1}}^{(g-1)-m}T_{{\bf n}_{m}}^{2(g-1)-m}\Big)
+\a_{X,\F_q}^{({\bf n}_{m-1})}(g-1)T_{{\bf n}_{m}}^{g-1})\Big)(1-T_{{\bf n}_{m}})(1-q_{{\bf n}_{m-1}}T_{{\bf n}_{m}})}{(1-T_{{\bf n}_{m}})(1-q_{{\bf n}_{m-1}}T_{{\bf n}_{m}})T_{{\bf n}_{m}}^{g-1}}\\
&\hskip 7.80cm+\frac{(q_{{\bf n}_{m-1}}-1)T_{{\bf n}_{m}}^{g}\b_{X,\F_q}^{({\bf n}_{m-1})}}{(1-T_{{\bf n}_{m}})(1-q_{{\bf n}_{m-1}}T_{{\bf n}_{m}})T_{{\bf n}_{m}}^{g-1}}\\
\ea
$$
and
$$
\Res_{T_{{\bf n}_{m}}=1}\wh\zeta_{X,\F_q}^{({\bf n}_{m-1})}(n_ms)=\Res_{T_{{\bf n}_{m}}=1}\frac{(q_{{\bf n}_{m-1}}-1)T_{{\bf n}_{m}}^{g}\b_{X,\F_q}^{({\bf n}_{m-1})}}{(1-T_{{\bf n}_{m}})(1-q_{{\bf n}_{m-1}}T_{{\bf n}_{m}})T_{{\bf n}_{m}}^{g-1}}
=\b_{X,\F_q}^{({\bf n}_{m-1})}=\wh v_{X,\F_q;1}^{({\bf n}_{m-1})}
$$
In addition, note that both 
$$\ba
&\wh\zeta_{X,\F_q}^{({\bf n}_{m-1})}(n_ms-n_m+1)\\
=&
\Big(\sum_{\ell=0}^{g-2}\a_{X,\F_q}^{({\bf n}_{m-1})}(\ell)\Big((q_{{\bf n}_{m-1}}^{n_m-1}T_{{\bf n}_{m}})^{m-(g-1)}+q_{{\bf n}_{m-1}}^{(g-1)-m}(q_{{\bf n}_{m-1}}^{n_m-1}T_{{\bf n}_{m}})^{(g-1)-m}\Big)\\
&+\a_{X,\F_q}^{({\bf n}_{m-1})}(g-1)(q_{{\bf n}_{m-1}}^{n_m-1}T_{{\bf n}_{m}})^{g-1})\Big)+\frac{(q-1)(q_{{\bf n}_{m-1}}^{n_m-1}T_{{\bf n}_{m}})^{g}\b_{X,\F_q}^{({\bf n}_{m-1})}}{(1-q_{{\bf n}_{m-1}}^{n_m-1}T_{{\bf n}_{m}})(1-q_{{\bf n}_{m-1}}q_{{\bf n}_{m-1}}^{n_m-1}T_{{\bf n}_{m}})(q_{{\bf n}_{m-1}}^{n_m-1}T_{{\bf n}_{m}})^{g-1}}\\
\ea
$$
and
$$
\frac{1}{(q_{{\bf n}_{m-1}}^{d_1+\ldots+d_{h-1}}T_{{\bf n}_{m}}-1)} 
$$
admit no pole at $T_{{\bf n}_{m}}=1$ unless $d_1+\ldots+d_{h-1}=0$ or equivalently, unless $d_h=n_m-1$, in which case 
$$
\frac{1}{(q_{{\bf n}_{m-1}}^{d_1+\ldots+d_{h-1}}T_{{\bf n}_{m}}-1)} =\frac{1}{T_{{\bf n}_{m}}-1}
$$ 
admits a simple pole at $T_{{\bf n}_{m}}=1$ with residue 1. Hence, by the functional equation,
$$
\Res_{T_{{\bf n}_{m}}=1}\frac{q_{{\bf n}_{m-1}}^{d_1+\ldots+d_{h-1}}T_{{\bf n}_{m}}\wh\zeta_{X,\F_q}^{({\bf n}_{m-1})}(n_ms-n_m+1)}
{(q_{{\bf n}_{m-1}}^{d_1+\ldots+d_{h-1}}T_{{\bf n}_{m}}-1)} =\bc 0&d_{h}<n_m-1\\
=\wh\zeta_{X,\F_q}^{({\bf n}_{m-1})}(n_m)&d_{h}=n_m-1\ec
$$
Finally, since both
$$
\ba
&\wh\zeta_{X,\F_q}^{({\bf n}_{m-1})}(n_ms-(d_1+\ldots+d_i))\\
=&
\Big(\sum_{\ell=0}^{g-2}\a_{X,\F_q}^{({\bf n}_{m-1})}(\ell)\Big((q_{{\bf n}_{m-1}}^{d_1+\ldots+d_i}T_{{\bf n}_{m}})^{\ell-(g-1)}+q_{{\bf n}_{m-1}}^{(g-1)-m}(q_{{\bf n}_{m-1}}^{d_1+\ldots+d_i}T_{{\bf n}_{m}})^{(g-1)-\ell}\Big)\\
+&\a_{X,\F_q}^{({\bf n}_{m-1})}(g-1)(q_{{\bf n}_{m-1}}^{d_1+\ldots+d_i}T_{{\bf n}_{m}})^{g-1})\Big)+\frac{(q_{{\bf n}_{m-1}}-1)(q_{{\bf n}_{m-1}}^{d_1+\ldots+d_i}T_{{\bf n}_{m}})^{g}\b_{X,\F_q}^{({\bf n}_{m-1})}}{(1-q_{{\bf n}_{m-1}}^{d_1+\ldots+d_i}T_{{\bf n}_{m}})(1-q_{{\bf n}_{m-1}}q_{{\bf n}_{m-1}}^{d_1+\ldots+d_i}T_{{\bf n}_{m}})(q_{{\bf n}_{m-1}}^{d_1+\ldots+d_i}T_{{\bf n}_{m}})^{g-1}}\\
\ea
$$
and
$$
\frac{1}{(q_{{\bf n}_{m-1}}^{d_1+\ldots+d_{i-1}}T_{{\bf n}_{m}}-1)(1-q_{{\bf n}_{m-1}}^{d_1+\ldots+d_{i+1}+1}T_{{\bf n}_{m}})}
$$
admit no pole at $T_{{\bf n}_{m}}=1$ unless 
$$
d_1+\ldots+d_{i-1}=0, 
$$ 
in which case,
$$
\frac{1}{(q_{{\bf n}_{m-1}}^{d_1+\ldots+d_{i-1}}T_{{\bf n}_{m}}-1)(1-q_{{\bf n}_{m-1}}^{d_1+\ldots+d_{i+1}+1}T_{{\bf n}_{m}})}=
\frac{1}{(T_{{\bf n}_{m}}-1)(1-q_{{\bf n}_{m-1}}^{d_i+d_{i+1}+1}T_{{\bf n}_{m}})}
$$
admits a simple pole with residue $\frac{1}{1-q_{{\bf n}_{m-1}}^{d_i+d_{i+1}+1}}$,
hence 
$$\ba
&\Res_{T_{{\bf n}_{m}}=1}\left(
\sum_{i=1}^{h-1}\frac{(1-q_{{\bf n}_{m-1}}^{d_i+d_{i+1}})\cdot q_{{\bf n}_{m-1}}^{d_1+\ldots+d_{i-1}}T_{{\bf n}_{m}}\wh\zeta_{X,\F_q}^{({\bf n}_{m-1})}(n_ms-(d_1+\ldots+d_i))}
{(q_{{\bf n}_{m-1}}^{d_1+\ldots+d_{i-1}}T_{{\bf n}_{m}}-1)(1-q_{{\bf n}_{m-1}}^{d_1+\ldots+d_{i+1}+1}T_{{\bf n}_{m}})}\right)\\
=&\bc 0&d_1+\ldots+d_{i-1}>0\\
\left(
\sum_{i=1}^{h-1}\frac{(1-q_{{\bf n}_{m-1}}^{d_i+d_{i+1}})\wh\zeta_{X,\F_q}^{({\bf n}_{m-1})}(d_i+1)}
{(1-q_{{\bf n}_{m-1}}^{d_i+d_{i+1}+1})}\right)&d_1+\ldots+d_{i-1}=0\ec\\
=&\bc 0&d_1+\ldots+d_{i-1}>0\\
\left(
\frac{(1-q_{{\bf n}_{m-1}}^{d_i+d_{i+1}})\wh\zeta_{X,\F_q}^{({\bf n}_{m-1})}(d_i+1)}
{(1-q_{{\bf n}_{m-1}}^{d_i+d_{i+1}+1})}\right)&d_1+\ldots+d_{i-1}=0\ec
\ea
$$
provided $i=1$. Therefore,
$$\ba
\b_{X,\F_q;{\bf n}_{m}}
=&q_{{\bf n}_{m-1}}^{\binom{n_m}{2}(g-1)}\Bigg(\sum_{\substack{d_1,\ldots,d_h>0\\ d_1+\ldots+d_h=n_m-1
}}\frac{\wh v^{({\bf n}_{m-1})}_{d_1}\ldots\wh v^{({\bf n}_{m-1})}_{d_h}}{\prod_{j=1}^{h-1}(1-q_{{\bf n}_{m-1}}^{d_j+d_{j+1}})}
\left(\Res_{T_{{\bf n}_{m}}=1}\frac{\wh\zeta_{X,\F_q}^{({\bf n}_{m-1})}(n_ms)}{(1-q_{{\bf n}_{m-1}}^{d_1+1}T_{{\bf n}_{m}})}\right.\\
&+\Res_{T_{{\bf n}_{m}}=1}\frac{q_{{\bf n}_{m-1}}^{d_1+\ldots+d_{h-1}}T_{{\bf n}_{m}}\wh\zeta_{X,\F_q}^{({\bf n}_{m-1})}(n_ms-n_m+1)}
{(q_{{\bf n}_{m-1}}^{d_1+\ldots+d_{h-1}}T_{{\bf n}_{m}}-1)}\\
&\left.+\Res_{T_{{\bf n}_{m}}=1}\left(
\sum_{i=1}^{h-1}\frac{(1-q_{{\bf n}_{m-1}}^{d_i+d_{i+1}})\cdot q_{{\bf n}_{m-1}}^{d_1+\ldots+d_{i-1}}T_{{\bf n}_{m}}\wh\zeta_{X,\F_q}^{({\bf n}_{m-1})}(n_ms-(d_1+\ldots+d_i))}
{(q_{{\bf n}_{m-1}}^{d_1+\ldots+d_{i-1}}T_{{\bf n}_{m}}-1)(1-q_{{\bf n}_{m-1}}^{d_1+\ldots+d_{i+1}+1}T_{{\bf n}_{m}})}\right)\right)\Bigg)\\
=&q_{{\bf n}_{m-1}}^{\binom{n_m}{2}(g-1)}\Bigg(\sum_{\substack{d_1,\ldots,d_h>0\\ d_1+\ldots+d_h=n_m-1
}}\frac{\wh v^{({\bf n}_{m-1})}_{d_1}\ldots\wh v^{({\bf n}_{m-1})}_{d_h}}{\prod_{j=1}^{h-1}(1-q_{{\bf n}_{m-1}}^{d_j+d_{j+1}})}
\frac{\wh v^{({\bf n}_{m-1})}_1}{(1-q_{{\bf n}_{m-1}}^{d_1+1})}+\wh\zeta_{X,\F_q}^{({\bf n}_{m-1})}(n)\wh v^{({\bf n}_{m-1})}_{n_m-1}\\
&\hskip 1.50cm+\sum_{\substack{d_1,\ldots,d_h>0\\ d_1+\ldots+d_h=n_m-1
}}\frac{\wh v^{({\bf n}_{m-1})}_{d_1}\ldots\wh v^{({\bf n}_{m-1})}_{d_h}}{\prod_{j=1}^{h-1}(1-q_{{\bf n}_{m-1}}^{d_j+d_{j+1}})}
\left(
\frac{(1-q_{{\bf n}_{m-1}}^{d_1+d_{2}})\wh\zeta_{X,\F_q}^{({\bf n}_{m-1})}(d_1+1)}
{(1-q_{{\bf n}_{m-1}}^{d_1+d_{2}+1})}\right)
\Bigg)\\
=&q_{{\bf n}_{m-1}}^{\binom{n_m}{2}(g-1)}\Bigg(\sum_{\substack{d_1,\ldots,d_h>0\\ d_1+\ldots+d_h=n_m-1
}}\frac{\wh v^{({\bf n}_{m-1})}_{d_1}\ldots\wh v^{({\bf n}_{m-1})}_{d_h}}{\prod_{j=1}^{h-1}(1-q_{{\bf n}_{m-1}}^{d_j+d_{j+1}})}
\frac{\wh v^{({\bf n}_{m-1})}_1}{(1-q_{{\bf n}_{m-1}}^{d_1+1})}\\
&\hskip 1.50cm+\wh v^{({\bf n}_{m-1})}_{n_m}+\sum_{d_1+1=a=2}^{n_m-1}\sum_{\substack{d_2,\ldots,d_h>0\\ d_1+\ldots+d_h=n_m-(d_1+1)
}}\frac{\wh v^{({\bf n}_{m-1})}_{d_2}\ldots\wh v^{({\bf n}_{m-1})}_{d_h}}{\prod_{j=2}^{h-1}(1-q_{{\bf n}_{m-1}}^{d_j+d_{j+1}})}
\left(\frac{\wh\zeta_{X,\F_q}^{({\bf n}_{m-1})}(d_1+1)}
{(1-q_{{\bf n}_{m-1}}^{d_{2}+(d_1+1)})}\right)
\Bigg)\\
=&q_{{\bf n}_{m-1}}^{\binom{n_m}{2}(g-1)}\Bigg(\sum_{\substack{d_1,\ldots,d_h>0\\ d_1+\ldots+d_h=n_m-1
}}\frac{\wh v^{({\bf n}_{m-1})}_{d_1}\ldots\wh v^{({\bf n}_{m-1})}_{d_h}}{\prod_{j=1}^{h-1}(1-q_{{\bf n}_{m-1}}^{d_j+d_{j+1}})}
\frac{\wh v^{({\bf n}_{m-1})}_1}{(1-q_{{\bf n}_{m-1}}^{d_1+1})}+\wh v^{({\bf n}_{m-1})}_{n_m}\\\ea$$
$$\ba
&\hskip 2.0cm+\sum_{a=2}^{n-1}\sum_{\substack{d_2,\ldots,d_h>0\\ d_2+\ldots+d_h=n_m-a
}}\frac{\wh v^{({\bf n}_{m-1})}_{d_2}\ldots\wh v^{({\bf n}_{m-1})}_{d_h}}{\prod_{j=2}^{h-1}(1-q_{{\bf n}_{m-1}}^{d_j+d_{j+1}})}
\left(\frac{\wh\zeta_{X,\F_q}^{({\bf n}_{m-1})}(a)\wh v^{({\bf n}_{m-1})}_{a-1}}
{(1-q_{{\bf n}_{m-1}}^{d_{2}+a})}\right)
\Bigg)\\
=&q_{{\bf n}_{m-1}}^{\binom{n_m}{2}(g-1)}\sum_{\substack{d_1,\ldots,d_h>0\\ d_1+\ldots+d_k=n_m
}}\frac{\wh v^{({\bf n}_{m-1})}_{d_1}\ldots\wh v^{({\bf n}_{m-1})}_{d_k}}{\prod_{j=1}^{k-1}(1-q_{{\bf n}_{m-1}}^{d_j+d_{j+1}})},\\
\ea$$
since, easily, we have
$$\ba
&\sum_{\substack{d_0,d_1,\ldots,d_k>0\\ d_0+d_1+\ldots+d_k=n_m
}}\frac{\wh v^{({\bf n}_{m-1})}_{d_0}\wh v^{({\bf n}_{m-1})}_{d_1}\ldots\wh v^{({\bf n}_{m-1})}_{d_k}}{\prod_{j=0}^{k-1}(1-q_{{\bf n}_{m-1}}^{d_j+d_{j+1}})}
=\sum_{d_0=1}^{n_m}\sum_{\substack{d_1,\ldots,d_k>0\\ d_1+\ldots+d_k=n_m-d_0
}}\frac{\wh v^{({\bf n}_{m-1})}_{d_1}\ldots\wh v^{({\bf n}_{m-1})}_{d_k}}{\prod_{j=1}^{k-1}(1-q_{{\bf n}_{m-1}}^{d_j+d_{j+1}})}\frac{\wh v^{({\bf n}_{m-1})}_{d_0}}{(1-q_{{\bf n}_{m-1}}^{d_0+d_1})}\\
=&\sum_{a=1}^{n_m}\wh v^{({\bf n}_{m-1})}_a\sum_{\substack{d_1,\ldots,d_k>0\\ d_1+\ldots+d_k=n_m-a
}}\frac{\wh v^{({\bf n}_{m-1})}_{d_1}\ldots\wh v^{({\bf n}_{m-1})}_{d_k}}{\prod_{j=1}^{k-1}(1-q_{{\bf n}_{m-1}}^{d_j+d_{j+1}})}\frac{1}{(1-q_{{\bf n}_{m-1}}^{a+d_1})}\\
=&\wh v^{({\bf n}_{m-1})}_{n_m}+\sum_{\substack{d_1,\ldots,d_k>0\\ d_1+\ldots+d_k=n_m-1
}}\frac{\wh v^{({\bf n}_{m-1})}_{n_1}\ldots\wh v^{({\bf n}_{m-1})}_{n_k}}{\prod_{j=1}^{k-1}(1-q_{{\bf n}_{m-1}}^{d_j+d_{j+1}})}\frac{\wh v^{({\bf n}_{m-1})}_1}{(1-q_{{\bf n}_{m-1}}^{1+n_1})}\\
&\hskip 2.0cm+\sum_{a=2}^{n-1}\wh v^{({\bf n}_{m-1})}_a\sum_{\substack{d_1,\ldots,d_k>0\\ d_1+\ldots+d_k=n_m-a
}}\frac{\wh v^{({\bf n}_{m-1})}_{n_1}\ldots\wh v^{({\bf n}_{m-1})}_{n_k}}{\prod_{j=1}^{k-1}(1-q_{{\bf n}_{m-1}}^{d_j+d_{j+1}})}\frac{1}{(1-q_{{\bf n}_{m-1}}^{a+n_1})}.\\
\ea$$
This then completes our proof of Theorem\,\ref{thm2.4}, and  the following closed formula for the ${\bf n}_m$-derived beta invariants $\b_{X,\F_q;{\bf n}_{m}}$ and hence the remaining part of Corollary\,\ref{cor2.5} about the beta invariant. 

\begin{thm}\label{thm2.6}
Let $X$ be an integral regular projective curve of genus $g$ over a finite field $\F_q$. Then
\be\label{eq46}
\b_{X,\F_q;{\bf n}_{m}}=q_{{\bf n}_{m-1}}^{\binom{n_m}{2}(g-1)}\sum_{\substack{d_1,\ldots,d_k>0\\ d_1+\ldots+d_k=n_m}}\frac{\wh v^{({\bf n}_{m-1})}_{n_1}\ldots\wh v^{({\bf n}_{m-1})}_{n_k}}{\prod_{j=1}^{k-1}(1-q_{{\bf n}_{m-1}}^{d_j+d_{j+1}})}.
\ee
\end{thm}

\subsection{Special counting miracle}\label{sec2.5}
By definition, the  ${\bf n}_m$-derived alpha and beta invariants 
$$
\a_{X,\F_q;{\bf n}_{m}}(\ell)\qquad(\ell=0,\ldots,g-1)\qqan \b_{X,\F_q;{\bf n}_{m}}
$$
determine the ${\bf n}_m$-derived zeta function $\wh\zeta_{X,\F_q}^{\,({\bf n}_m)}(s)$ uniquely.
In this subsection, we expose an important relation among them. For this, we set
$$
{\bf n}_m+1:=(n_0,n_1,\ldots, n_m+1)\quad\forall {\bf n}_m:=(n_0,n_1,\ldots, n_m).
$$
\begin{thm}[Special Counting Miracle]\label{thm2.7} Let $X$ be an integral regular projective curve of genus $g$. Then
\be
\a_{X,\F_q}^{({\bf n}_{m}+1)}(0)=q_{{\bf n}_{m-1}}^{n^{~}_m(g-1)}\a_{X,\F_q}^{({\bf n}_{m-1})}(0)\cdot\b_{X,\F_q;{\bf n}_{m}}
\ee
\end{thm}
\begin{rmk} In the classical setting, namely $m=0$, this special counting miracle  is conjectured by the author in \cite{W3}, and proved first for elliptic curves by Zagier and myself in \cite{WZ1}, and later established by Sugahara (see the appendix in \cite{WRH}). An independent proof for the classical can be found in  \cite{MR} and \cite{WRH} as well.
\end{rmk}
\bp
By Corollary\,\ref{cor2.5},
$\a_{X,\F_q;{\bf n}_{m}}(0)$ is the  constant term of the polynomial
$$
\Big(1-T_{{\bf n}_{m}}\Big)\Big(1-q_{{\bf n}_{m}}T_{{\bf n}_{m}}\Big)\cdot T_{{\bf n}_{m}}^{g-1}\wh Z_{X,\F_q}^{\,({\bf n}_{m})}(T_{{\bf n}_{m}})
$$
This latest polynomial, by Definition\,\ref{defn1.1} is given
$$\ba
&\Big(1-T_{{\bf n}_{m}}\Big)\Big(1-q_{{\bf n}_{m}}T_{{\bf n}_{m}}\Big)\\
&\times q_{{\bf n}_{m-1}}^{\binom{n^{~}_m}{2}(g-1)}\sum_{a=1}^{n_m}\sum_{\substack{k_1,\ldots,k_p>0\\ k_1+\ldots+k_p=n_m-a}}\frac{\wh v_{X,\F_q,k_1}^{({\bf n}_{m-1})}\ldots\wh v_{X,\F_q,k_p}^{({\bf n}_{m-1})}}{\prod_{j=1}^{p-1}(1-q_{{\bf n}_{m-1}}^{k_j+k_{j+1}})} \frac{1}{(1-q_{{\bf n}_{m-1}}^{-n_m+a+k_{p}}T_{{\bf n}_{m}}^{-1})}\\
 &\hskip 0.20cm\times T_{{\bf n}_{m}}^{g-1}\wh\zeta_{X,\F_q}^{\,({\bf n}_{m-1})}(n_ms-n_m+a)\sum_{\substack{l_1,\ldots,l_r>0\\ l_1+\ldots+l_r=a-1}}
 \frac{1}{(1-q_{{\bf n}_{m-1}}^{n_m-a+1+l_{1}}T_{{\bf n}_{m}})}\frac{\wh v_{X,\F_q,l_1}^{({\bf n}_{m-1})}\ldots\wh v_{X,\F_q,l_r}^{({\bf n}_{m-1})}}{\prod_{j=1}^{r-1}(1-q_{{\bf n}_{m-1}}^{l_j+l_{j+1}})}\\
\ea$$
By examining the summand in the summation
$\sum_{\substack{k_1,\ldots,k_p>0\\ k_1+\ldots+k_p=n_m-a}}$, particularly, the factor
$$
\frac{1}{(1-q_{{\bf n}_{m-1}}^{-n_m+a+k_{p}}T_{{\bf n}_{m}}^{-1})}=\frac{T_{{\bf n}_{m}}}{(T_{{\bf n}_{m}}-q_{{\bf n}_{m-1}}^{-n_m+a+k_{p}})}
$$
we conclude that the non-zero contribution from the summands in the summation $\sum_{a=1}^{n_m}$ comes only from $a=n_m$. So we only need to value of the following rational function at $0$:
$$\ba
&\Big(1-T_{{\bf n}_{m}}\Big)\Big(1-q_{{\bf n}_{m}}T_{{\bf n}_{m}}\Big)\\
&\times q_{{\bf n}_{m-1}}^{\binom{n^{~}_m}{2}(g-1)} T_{{\bf n}_{m}}^{g-1}\wh\zeta_{X,\F_q}^{\,({\bf n}_{m-1})}(n_ms)\sum_{\substack{l_1,\ldots,l_r>0\\ l_1+\ldots+l_r=n_m-1}}
 \frac{1}{(1-q_{{\bf n}_{m-1}}^{1+l_{1}}T_{{\bf n}_{m}})}\frac{\wh v_{X,\F_q,l_1}^{({\bf n}_{m-1})}\ldots\wh v_{X,\F_q,l_r}^{({\bf n}_{m-1})}}{\prod_{j=1}^{r-1}(1-q_{{\bf n}_{m-1}}^{l_j+l_{j+1}})}\\
\ea$$
Recall that
$$\hskip -1.5cm\ba
T_{{\bf n}_{m}}^{g-1}\wh\zeta_{X,\F_q}^{\,({\bf n}_{m-1})}(n_ms)
=&\Big(\sum_{\ell=0}^{g-2}\a_{X,\F_q}^{({\bf n}_{m-1})}(\ell)\Big(T_{{\bf n}_{m}}^{\ell}+q_{{\bf n}_{m}}^{(g-1)-\ell}T_{{\bf n}_{m}}^{2(g-1)-\ell}\Big)
+\a_{X,\F_q}^{({\bf n}_{m-1})}(g-1)T_{{\bf n}_{m}}^{2(g-1)})\Big)\\
&\hskip 5.0cm+\frac{(q_{{\bf n}_{m}}-1)T_{{\bf n}_{m}}^{g}\b_{X,\F_q;{\bf n}_{m-1}}}{(1-T_{{\bf n}_{m}})(1-q_{{\bf n}_{m}}T_{{\bf n}_{m}})}
\ea$$
Therefore,
$$\ba
\a_{X,\F_q;{\bf n}_{m}}(0)
=&q_{{\bf n}_{m-1}}^{\binom{n^{~}_m}{2}(g-1)}\a_{X,\F_q}^{({\bf n}_{m-1})}(0)
\sum_{\substack{l_1,\ldots,l_r>0\\ l_1+\ldots+l_r=n_m-1}}
\frac{\wh v_{X,\F_q,l_1}^{({\bf n}_{m-1})}\ldots\wh v_{X,\F_q,l_r}^{({\bf n}_{m-1})}}{\prod_{j=1}^{r-1}(1-q_{{\bf n}_{m-1}}^{l_j+l_{j+1}})}\\
=&q_{{\bf n}_{m-1}}^{\binom{n^{~}_m}{2}(g-1)}\a_{X,\F_q}^{({\bf n}_{m-1})}(0)\b_{X,\F_q;{\bf n}_{m}-1}''
=q_{{\bf n}_{m-1}}^{(n^{~}_m-1)(g-1)}\a_{X,\F_q}^{({\bf n}_{m-1})}(0)\b_{X,\F_q;{\bf n}_{m-1},n_m-1}\\
=&q_{{\bf n}_{m-1}}^{(n^{~}_m-1)(g-1)}\a_{X,\F_q}^{({\bf n}_{m-1})}(0)\b_{X,\F_q;{\bf n}_{m}-1}
\ea$$
as wanted.
\ep

As commented at the end of \S\ref{sec1.3}, the factor $\a_{X,\F_q}^{({\bf n}_{m-1})}(0)$ actually introduce an additional $\a_{X,\F_q}^{({\bf n}_{m-1})}(0)^{n_m}$ factor to $\a_{X,\F_q;{\bf n}_{m}}(0)$. This in effect says that there is an unessential factor $\a_{X,\F_q}^{({\bf n}_{m-1})}(0)^{n_m-1}$ involved in  $\b_{X,\F_q;{\bf n}_{m}-1}$.

\subsection{Positivity conjecture and its application}\label{sec2.6}
The ${\bf n}_m$-derived alpha and beta invariants $\a_{X,\F_q;{\bf n}_{m}}(\ell),\ \b_{X,\F_q;{\bf n}_{m}}$ plays a key role  in the theory of the  ${\bf n}_m$-derived zeta functions $\wh\zeta_{X,\F_q}^{\,({\bf n}_{m})}(s)$ for curves over finite fields. 
In this subsection, we formulate the following
\begin{conj}[Positivity of ${\bf n}_m$-derived alpha and beta invariants]\label{conj2.1} Let $X$ be an integral regular projective curve of genus $g$ over $\F_q$. Then the
${\bf n}_m$-derived alpha and beta invariants 
\be
\a_{X,\F_q;{\bf n}_{m}}(\ell)\qquad(\ell=0,\ldots,g-1)\qqan \b_{X,\F_q;{\bf n}_{m}}
\ee are always strictly positive.
\end{conj}
The positivity in classical situation is rather trivial since we have the following geometric interpretation of these non-abelian invariants:
\be
\a_{X,\F_q;n}(d)=\sum_{\cE}\frac{q^{h^0(X,\cE)}-1}{\#\Aut(\cE)}\qqan \b_{X,\F_q;n}(d)=\sum_{\cE}\frac{1}{\#\Aut(\cE)}
\ee
where $\cE$ runs over all rank $n$ semi-stable $\F_q$-rational vector bundle of degree $d$. 

Recall that, by Definition\,\ref{defn2.1}, for any $1\leq n\leq n_m$, the  rational function
$\De_{X,\F_q;{\bf n}_{m};n}(T_{{\bf n}_m})$,  resp. the polynomial $\Ga_{X,\F_q;{\bf n}_{m};n}(T_{{\bf n}_m})$,  of $T_{{\bf n}_m}$ is defined by
$$\ba
\De_{X,\F_q;{\bf n}_{m};n}(T_{{\bf n}_m}):=&\sum_{\substack{k_1,\ldots,k_p>0\\ k_1+\ldots+k_p=n}}\frac{\wh v_{X,\F_q,k_1}^{({\bf n}_{m-1})}\ldots\wh v_{X,\F_q,k_p}^{({\bf n}_{m-1})}}{\prod_{j=1}^{p-1}\Big(1-q_{{\bf n}_{m-1}}^{k_j+k_{j+1}}\Big)} \frac{1}{\Big(q_{{\bf n}_{m-1}}^{-n_ms+k_{p}}-1\Big)}\\
=&
\sum_{\substack{k_1,\ldots,k_p>0\\ k_1+\ldots+k_p=n}}\frac{\wh v_{X,\F_q,k_1}^{({\bf n}_{m-1})}\ldots\wh v_{X,\F_q,k_p}^{({\bf n}_{m-1})}}{\prod_{j=1}^{p-1}\Big(1-q_{{\bf n}_{m-1}}^{k_j+k_{j+1}}\Big)} \frac{1}{\Big(q_{{\bf n}_{m-1}}^{k_{p}}T_{{\bf n}_{m}}-1\Big)}
\ea
$$
resp.
$$
\Ga_{X,\F_q;{\bf n}_{m},n}(T_{{\bf n}_m}):=\De_{X,\F_q;{\bf n}_{m-1},n}(T_{{\bf n}_m})\cdot\prod_{\substack{\ell=1}}^n\Big(q_{{\bf n}_{m-1}}^{\ell} T_{{\bf n}_{m}}-1\Big).
$$
Our main result of  this section is the following
\begin{thm}\label{thm5.1} All the roots of the polynomial $\Ga_{X,\F_q;{\bf n}_{m},n}(T_{{\bf n}_m})$ are real. Furthermore, there exists one and only one root of $\Ga_{X,\F_q;{\bf n}_{m},n}(T_{{\bf n}_m})$ in each open interval 
\be
\Big(q_{{\bf n}_{m-1}}^{-\kappa-1}, q_{{\bf n}_{m-1}}^{-\kappa}\Big)\qquad\kappa=1,2,\ldots,n-1
\ee 
provided that the Positivity Conjecture at the level $(m-1)$, i.e., Conjecture\,\ref{conj2.1} on positivity of the ${\bf n}_{m-1}$-derived alpha and beta invariants, holds. 
\end{thm}
\bp
We begin with the following 
\begin{lem}\label{lem5.2} Let $k$ be a (strictly) positive integer. Then
$\wh\zeta^{\,({\bf n}_{m-1})}(k)$  and hence $\wh v_{X,\F_q,k}^{({\bf n}_{m-1})}$ are (strictly) positive, provided that the Positivity Conjecture holds at the level $(m-1)$.
\end{lem} 
\bp
By Corollary\,\ref{cor2.5},
$$\ba
&\wh \zeta_{X,\F_q}^{({\bf n}_{m-1})}(k)=\wh \zeta_{X,\F_q}^{({\bf n}_{m-1})}(1-k)=\frac{(q_{{\bf n}_{m-1}}-1)q_{{\bf n}_{m-1}}^{(k-1)}\b_{X,\F_q;{\bf n}_{m-1}}}{(q^{k-1}_{{\bf n}_{m-1}}-1)(q_{{\bf n}_{m-1}}q^{k-1}_{{\bf n}_{m-1}}-1)}\\
&+\Big(\sum_{\ell=0}^{g-2}\a_{X,\F_q}^{({\bf n}_{m-1})}(\ell)\Big(q_{{\bf n}_{m-1}}^{(k-1)(\ell-(g-1))}+q_{{\bf n}_{m-1}}^{(g-1)-\ell}q_{{\bf n}_{m-1}}^{(k-1)((g-1)-\ell)}\Big)
+\a_{X,\F_q}^{({\bf n}_{m-1})}(g-1)q_{{\bf n}_{m-1}}^{(k-1)(g-1)})\Big)\\
\ea$$
Hence $\wh \zeta_{X,\F_q}^{({\bf n}_{m-1})}(k)$ is strictly positive when $k\geq 2$ since each term involved is strictly positive by the level $(m-1)$ positivity conjecture.
Similarly,
\be
\wh \zeta_{X,\F_q}^{({\bf n}_{m-1})}(1)=\Res_{T_{{\bf n}_{m-1}}=1}\wh Z_{X,\F_q}^{\,({\bf n}_{m-1})}(T_{{\bf n}_{m-1}})=\b_{X,\F_q;{\bf n}_{m-1}}
\ee
which  is strictly positive by the level $(m-1)$ positivity conjecture.
\ep
Indeed, by definition, we have, 
$$\ba
\Ga_{X,\F_q;{\bf n}_{m},n}(T_{{\bf n}_m})
=&
\sum_{\substack{k_1,\ldots,k_p>0\\ k_1+\ldots+k_p=n}}\frac{\wh v_{X,\F_q,k_1}^{({\bf n}_{m-1})}\ldots\wh v_{X,\F_q,k_p}^{({\bf n}_{m-1})}}{\prod_{j=1}^{p-1}\Big(1-q_{{\bf n}_{m-1}}^{k_j+k_{j+1}}\Big)} \frac{1}{\Big(q_{{\bf n}_{m-1}}^{k_{p}}T_{{\bf n}_{m}}-1\Big)}\cdot\prod_{\substack{\ell=1}}^n\Big(q_{{\bf n}_{m-1}}^{\ell} T_{{\bf n}_{m}}-1\Big)\\
=&
\sum_{\substack{k_1,\ldots,k_p>0\\ k_1+\ldots+k_p=n}}\frac{\wh v_{X,\F_q,k_1}^{({\bf n}_{m-1})}\ldots\wh v_{X,\F_q,k_p}^{({\bf n}_{m-1})}}{\prod_{j=1}^{p-1}\Big(1-q_{{\bf n}_{m-1}}^{k_j+k_{j+1}}\Big)}\cdot\prod_{\substack{\ell=0,\ell\not=k_p}}^n\Big(q_{{\bf n}_{m-1}}^{\ell} T_{{\bf n}_{m}}-1\Big)\\
\ea
$$
This is a degree $n-1$ polynomial in $T_{{\bf n}_m}$ with real coefficients.
Now, for $\kappa=1,2,\ldots n$,
$$\ba
&\Ga_{X,\F_q;{\bf n}_{m},n}(q_{{\bf n}_{m-1}}^{-\kappa})
=
\sum_{\substack{k_1,\ldots,k_{p-1}>0, k_P=\kappa\\ k_1+\ldots+k_p=n}}\frac{\wh v_{X,\F_q,k_1}^{({\bf n}_{m-1})}\ldots\wh v_{X,\F_q,k_p}^{({\bf n}_{m-1})}}{\prod_{j=1}^{p-1}\Big(1-q_{{\bf n}_{m-1}}^{k_j+k_{j+1}}\Big)}\cdot\prod_{\substack{\ell=1,\ell\not=k_p=\kappa}}^n\Big(q_{{\bf n}_{m-1}}^{\ell-\kappa} -1\Big)\\
=&
\sum_{\substack{k_1,\ldots,k_{p-1}>0, k_P=\kappa\\ k_1+\ldots+k_{p-1}=n-\kappa}}\frac{\wh v_{X,\F_q,k_1}^{({\bf n}_{m-1})}\ldots\wh v_{X,\F_q,k_{p-1}}^{({\bf n}_{m-1})}}{\prod_{j=1}^{p-2}\Big(1-q_{{\bf n}_{m-1}}^{k_j+k_{j+1}}\Big)}\frac{\wh v_{X,\F_q,\kappa}^{({\bf n}_{m-1})}}{1-q_{{\bf n}_{m-1}}^{k_{p-1}+\kappa}}\cdot\prod_{\substack{\ell=1,\ell\not=k_p=\kappa}}^n\Big(q_{{\bf n}_{m-1}}^{\ell-\kappa} -1\Big)\\
=&(-1)^{\kappa-1}\wh v_{X,\F_q,\kappa}^{({\bf n}_{m-1})}
\sum_{\substack{k_1,\ldots,k_{p-1}>0, k_P=\kappa\\ k_1+\ldots+k_{p-1}=n-\kappa}}\frac{\wh v_{X,\F_q,k_1}^{({\bf n}_{m-1})}\ldots\wh v_{X,\F_q,k_{p-1}}^{({\bf n}_{m-1})}}{\prod_{j=1}^{p-2}\Big(1-q_{{\bf n}_{m-1}}^{k_j+k_{j+1}}\Big)}\frac{1}{1-q_{{\bf n}_{m-1}}^{k_{p-1}+\kappa}}\cdot\prod_{\ell=1}^{\kappa-1}\Big(1-q_{{\bf n}_{m-1}}^{\ell-\kappa} \Big)\prod_{\ell=\kappa+1}^n\Big(q_{{\bf n}_{m-1}}^{\ell-\kappa} -1\Big)\\
=&(-1)^{\kappa}\wh v_{X,\F_q,\kappa}^{({\bf n}_{m-1})}\cdot \Ga_{X,\F_q;{\bf n}_{m},n-\kappa}(q_{{\bf n}_m}^{\kappa})
\cdot\prod_{\ell=1}^{\kappa-1}\Big(1-q_{{\bf n}_{m-1}}^{\ell-\kappa} \Big)\prod_{\ell=\kappa+1}^n\Big(q_{{\bf n}_{m-1}}^{\ell-\kappa} -1\Big)\\
\ea$$
We claim that this sign is simply $(-1)^{\kappa}$. Indeed, with this latest relation, an induction in $n$ can be applied to conclude that the sign of $\Ga_{X,\F_q;{\bf n}_{m},n}(q_{{\bf n}_{m-1}}^{-\kappa})$ is given by $(-1)^{\kappa+1}$. To see this, first we note that
$$\ba
\Ga_{X,\F_q;{\bf n}_{m},n}(0)=&
\sum_{\substack{k_1,\ldots,k_p>0\\ k_1+\ldots+k_p=n}}\frac{\wh v_{X,\F_q,k_1}^{({\bf n}_{m-1})}\ldots\wh v_{X,\F_q,k_p}^{({\bf n}_{m-1})}}{\prod_{j=1}^{p-1}\Big(1-q_{{\bf n}_{m-1}}^{k_j+k_{j+1}}\Big)}\cdot\prod_{\substack{\ell=1,\ell\not=k_p}}^n\Big(q_{{\bf n}_{m-1}}^{\ell} 0-1\Big)\\
=&(-1)^{n-1}\sum_{\substack{k_1,\ldots,k_p>0\\ k_1+\ldots+k_p=n}}\frac{\wh v_{X,\F_q,k_1}^{({\bf n}_{m-1})}\ldots\wh v_{X,\F_q,k_p}^{({\bf n}_{m-1})}}{\prod_{j=1}^{p-1}\Big(1-q_{{\bf n}_{m-1}}^{k_j+k_{j+1}}\Big)}
\ea
$$
Thus for $1\leq a\leq n$,
the sign of $\Ga_{X,\F_q;{\bf n}_{m},a}(0)$, which is the same as the sign of $\Ga_{X,\F_q;{\bf n}_{m},a}(q_{{\bf n}_{m-1}}^{-a})$ by the inductive hypothesis,
 is $(-1)^{a+1}$. Consequently, by the inductive hypothesis again, we conclude that the sign of $\Ga_{X,\F_q;{\bf n}_{m},a}(q_{{\bf n}_{m-1}}^{-1})$ is $(-1)^{a+1}(-1)^{a-1}=1$, since there are sign changes at $q_{{\bf n}_{m-1}}^{-\kappa}$ for each $\kappa=1,2,\ldots, a$. Furthermore, since all $a-1$ roots of  $\Ga_{X,\F_q;{\bf n}_{m},a}(T_{{\bf n}_m})$ are located in the interval $(q_{{\bf n}_{m-1}}^{-a},q_{{\bf n}_{m-1}}^{\kappa})$, then sign of $\Ga_{X,\F_q;{\bf n}_{m},a}(q_{{\bf n}_{m-1}}^{-1})$ is the same as that of $\Ga_{X,\F_q;{\bf n}_{m},a}(q_{{\bf n}_{m-1}}^{-1})$, namely, $1$.
Therefore, $\Ga_{X,\F_q;{\bf n}_{m},n-\kappa}(q_{{\bf n}_m}^{\kappa})
$ is positive, since by Lemma\,\ref{lem5.2}, $\wh v_{X,\F_q,\kappa}^{({\bf n}_{m-1})}$ is (strictly) positive, then the sign of $\Ga_{X,\F_q;{\bf n}_{m},n}(q_{{\bf n}_{m-1}}^{-\kappa})$ is $(-1)^{\kappa}$ as claimed above. Therefore there are changes of $\Ga_{X,\F_q;{\bf n}_{m},n}(T_{{\bf n}_m})$ at $q_{{\bf n}_{m-1}}^{-\kappa}$ for $\kappa=1,2,\ldots,n$. This then completes the proof since $\Ga_{X,\F_q;{\bf n}_{m},n}(T_{{\bf n}_m})$  is a degree $n-1$-polynomial in $T_{{\bf n}_m}$.
\ep

\subsection{${\bf n}_m$-derived Riemann hypothesis}

The most surprising point for the ${\bf n}_m$-derived zeta functions of curves over finite field is that we have the following
\begin{conj} [${\bf n}_m$-Derived Riemann Hypothesis]\label{conjA} Let $X$ be an integral regular projective curve of genus $g$ over $\F_q$. Then the
${\bf n}_m$-derived zeta function $\wh \zeta^{\,({\bf n}_m)}_{X,\F_q}(s)$ satisfies the Riemann Hypothesis. That is, all zeros of   $\wh \zeta^{\,({\bf n}_m)}_{X,\F_q}(s)$ lie on the central line $\Re(s)=\frac{1}{2}$.
\end{conj}
Obviously, this ${\bf n}_m$-derived Riemann hypothesis is equivalent to the fact that
all reciprocity roots $\a_{X,\F_q;\ell}^{\,({\bf n}_m)}$ $(1\leq \ell\leq 2g)$ of ${\bf n}_m$-derived Zeta function $\wh Z^{\,({\bf n}_m)}(s)$ satisfies the condition that
\be
\left|\a_{X,\F_q;\ell}^{\,({\bf n}_m)}\right|=q_{{\bf n}_m}^{\frac{1}{2}}\qquad (1\leq \ell\leq 2g).
\ee
We will prove in the final section that the  following
\begin{thm} Let $X$ be an integral regular projective curve of genus $g$ over $\F_q$. Then the
${\bf n}_m$-derived Riemann hypothesis holds for  the
${\bf n}_m$-derived zeta function $\wh \zeta_{X,\F_q}^{\,({\bf n}_m)}(s)$, for any of the following cases:

\begin{enumerate}
\item[(1)] $g=1$. That is, $X$ is an elliptic curve $E$ over finite field. $\qquad$
(2) ${\bf n_m}=(2,2,\ldots,2)$ 
\end{enumerate}
\end{thm}

This indicates that,  being a crucial zeta property, Riemann hypothesis is indeed quite universal.
It admits vast diversities, despite the fact that the pioneers and later fundamentalists believed that the Riemann hypothesis is coming from the multiplicative structure of Euler products.
 
\section{Multiplicative Structure of ${\bf n}_m$-Derived Zeta Functions}
Before we go further, let us expose the hidden multiplicative structure for our ${\bf n}_m$-Derived zeta Functions of curves.
\subsection{Multiplicative structure of Artin zeta functions}
We start with $m=0$ and $n_0=1$. Then the associated 0-th (1)-derived zeta function for an integral regular projective curve $X$ of genus $g$ over $\F_q$ is simply the complete Artin zeta function $\zeta_{X,\F_q}(s)$ of $X$ over $\F_q$. As such,  we have
$$
q^{-s(g-1)}\wh\zeta_{X,\F_q}(s)=\zeta_{X,\F_q}(s):=\sum_{D\geq 0}N(D)^{-s}=\prod_{P\in X}\frac{1}{1-N(P)^{-s}}.
$$
where $D$ runs over all effective divisors on $X$, and $P$ runs through all closed points on $X$.
Furthermore, by applying the rationality and the functional equationwe get, with $t=q^{-s}$,
$$\ba
\zeta_{X,\F_q}(s)=&\frac{\prod_{\ell=1}^g(1-\a_{X,\F_q;\ell}t)(1- \b_{X,\F_q;\ell}t)}{(1-t)(1-qT)}\\
=&\exp\left(\sum_{\ell=1}^g\log(1-\a_{X,\F_q;\ell}t)+\log(1-\b_{X,\F_q;\ell}t)-\log(1-t)-\log(1-qt)\right)\\
=&\exp\left(\sum_{k=1}^\infty N_k\frac{t^k}{k}\right)\\
\ea
$$
where 
$$
N_k=q^k+1-\sum_{\ell=1}^g\Big(\a_{X,\F_q;\ell}^k+\b_{X,\F_q;\ell}^k\Big)
$$
From this expression it is not difficult to arrive at the following
\begin{lem} The Hasse-Weil theorem on the rank one Riemann hypothesis is equivalent to the condition that
$$
|N_k-q^k-1|\leq C\sqrt {q^k}\qquad (\forall k\gg 0)
$$
\end{lem}
In addition, it is not difficult to see that
$$
N_1=\#X(\F_q).
$$
Furthermore, by applying the above multiplicative structure in terms of Euler product (for details, see the next subsection), we conclude also that
$$
N_k=\#X(F_{q^k}).
$$
In fact, this relation is equivalent to the Euler product structure for Artin zeta functions.

\subsection{Multiplicative structure of ${\bf n}_m$-derived zeta functions}

It is a natural question whether the multiplicative structure exists for the rank $n$ non-abelian zeta functions of curves over finite fields. Our first level answer is no in the sense that
there is no possible to obtain an Euler product structure for rank $n$ zeta functions when $n\geq 2$, since Euler product is  commutative. However, this does
prevent us to go further to give a more refined studies on the multiplicative structure of the rank $n$, or more generally, of the ${\bf n}_m$-derived zeta functions for curves over finite fields. 
To put it plainly, we claim that 
\begin{enumerate}
\item [(1)] {\it there are natural multiplicative structures for every ${\bf n}_m$-derived zeta functions of  curves over finite fields} and 
\item[(2)] {\it this multiplicative structure plays a key role in our proof on the ${\bf n}_m$-derived Riemann Hypothesis on  ${\bf n}_m$-derived zeta functions for curves over finite fields}. 
\end{enumerate}

Indeed, by Theorems\,\ref{thm2.1},\,\ref{thm2.2},\,\ref{thm2.4}, we have
$$\ba
Z_{X,\F_q}^{({\bf n}_m)}(T_{{\bf n}_m}):=&T_{{\bf n}_m}^{(g-1)}\wh Z_{X,\F_q}^{\,({\bf n}_m)}(T_{{\bf n}_m})\\
=&L_{X,\F_q}^{({\bf n}_m)}\frac{\prod_{\ell=1}^g(1-\a_{X,\F_q;{\bf n}_m,\ell}T_{{\bf n}_m})(1-\b_{X,\F_q;{\bf n}_m,\ell}T_{{\bf n}_m})}{(1-T_{{\bf n}_m})(1-q_{{\bf n}_m}T_{{\bf n}_m})}\\
=&L_{X,\F_q}^{({\bf n}_m)}\exp\left(\sum_{k=1}^\infty N_{X,\F_q;k}^{({\bf n}_m)}\frac{T_{{\bf n}_m}^k}{k}\right)
\ea
$$
where
$$
N_{X,\F_q;k}^{({\bf n}_m)}=q_{{\bf n}_m}^k+1-\sum_{\ell=1}^g\Big(\a_{X,\F_q;{\bf n}_m,\ell}^k+\a_{X,\F_q;{\bf n}_m,\ell}^k\Big).
$$
In particular,
$$\ba
 L_{X,\F_q}^{({\bf n}_m)}\frac{N_{X,\F_q;1}^{({\bf n}_m)}}{q_{{\bf n}_m}-1}=&L_{X,\F_q}^{({\bf n}_m)}\frac{\prod_{\ell=1}^g(1-\a_{X,\F_q;{\bf n}_m,\ell})(1-\b_{X,\F_q;{\bf n}_m,\ell})}{q_{{\bf n}_m}-1}\\
=&\Res_{T_{{\bf n}_m}=1}Z_{X,\F_q}^{({\bf n}_m)}(T_{{\bf n}_m})=\Res_{T_{{\bf n}_m}=1}\wh Z_{X,\F_q}^{\,({\bf n}_m)}(T_{{\bf n}_m})\\
=&\b_{X,\F_q;{\bf n}_m}.
\ea
$$
Hence, if we let $\zeta_l$ is an $l$-th primitive root of unity, then
$$\ba
\prod_{k=1}^l \wh Z_{X,\F_q}^{({\bf n}_m)}(\zeta_l^k T_{{\bf n}_m}):=&T_{{\bf n}_m}^{(g-1)l}\prod_{k=1}^l \wh Z_{X,\F_q}^{\,({\bf n}_m)}(\zeta_l^kT_{{\bf n}_m})\\
=&(L_{X,\F_q}^{({\bf n}_m)})^l\frac{\prod_{\ell=1}^g(1-\a_{X,\F_q;{\bf n}_m,\ell}^lT_{{\bf n}_m}^l)(1-\b_{X,\F_q;{\bf n}_m,\ell}^lT_{{\bf n}_m}^l)}{(1-T_{{\bf n}_m}^l)(1-q_{{\bf n}_m}^lT_{{\bf n}_m}^l)}\\
=&(L_{X,\F_q}^{({\bf n}_m)})^l\exp\left(\sum_{k=1}^\infty N_{X,\F_q;kl}^{({\bf n}_m)}\frac{T_{{\bf n}_m}^{kl}}{k}\right)
\ea
$$
since
\be\label{eq60}
\prod_{k=1}^l(1-\zeta_l^kz)=1-z^l\qqan\sum_{k=1}^l(\zeta_l^m)^k=\bc l&l|m\\
0&l\not|m\ec
\ee
For the classical Artin zeta function, this then leads to the famous relation that
$$
Z_{X_l/\F_{q^l}}(t^l)=\prod_{k=1}^lZ_{X/\F_{q}}(\zeta_l^kt)
$$
with the help of the Euler product and the relation
$$
N_1(X_l)=N_l(X).
$$

However, except for the relation above \eqref{eq60}, for derived zetas, there is no clearly statement similar to that of the Artin zetas. Thus the best approximation of the multiplication structure for the ${\bf n}_m$-derived zeta function of curves over finite fields is based on the inductive definition, the ${\bf n}_m$-derived zeta functions are build up based on ${\bf n}_{m-1}$-derived zeta functions. In this way, we are lead to the $0$-th derived zeta functions whcih are nothing but the rank $n$ non-abelian zeta function $\wh \zeta_{X,\F_q;n}(s)$. Therefore, finally using the special uniformity of zetas, we obtain a multiplication structure from that of Artin zetas.

As such, the hidden multiplication structure of ${\bf n}_m$-derived zeta functions of curves over finite fields are quite complicated and implicit. However, another form of multiplicative structure ${\bf n}_m$-derived zeta functions of curves over finite fields exists based on the invariants $N^{({\bf n}_m)}_{X,\F_q;l}$. 

\subsection{Another multiplicative structure for ${\bf n}_m$-derived zeta functions}

Our  next aim is to expose this new form of multiplicative structure for ${\bf n}_m$-derived zeta functions of curves over finite fields, motivated by \cite{WZ1}.

To start with,  introduce then the well-defined infinite product
\be\label{eq64}
B_{X,\F_q}^{({\bf n}_m)}(x):=\exp\left(\sum_{m=1}^{\infty}\frac{N_{X,\F_q;m}^{({\bf n}_m)}}{q_{{\bf n}_m}^m-1}\frac{x^m}{m}\right)=\sum_{k=0}^\infty b_{X,\F_q;k}^{({\bf n}_m)}x^k.
\ee
Then
$$
b_{X,\F_q;k}^{({\bf n}_m)}=\sum_{\substack{k_1,k_2,\ldots\geq 1\\ k_1+2k_2+\cdots=k}}\frac{(N_{X,\F_q;1}^{({\bf n}_m)})^{k_1}(N_{X,\F_q;2}^{({\bf n}_m)})^{k_2}\cdots}{(q_{{\bf n}_m}-1)^{k_1}(2(q_{{\bf n}_m}^2-1))^{k_1}\cdots k_1!k_2!\cdots}
$$
and
$$\ba
\frac{B_{X,\F_q}^{({\bf n}_m)}(q_{{\bf n}_m}x)}{B_{X,\F_q}^{({\bf n}_m)}(x)}
=&\exp\left(\sum_{m=1}^{\infty}(q_{{\bf n}_m}^m-1)\frac{N_{X,\F_q;m}^{({\bf n}_m)}}{q_{{\bf n}_m}^m-1}\frac{x^m}{m}\right)
=\exp\left(\sum_{m=1}^{\infty}N_{X,\F_q;m}^{({\bf n}_m)}\frac{x^m}{m}\right)\\
=&\frac{\prod_{\ell=1}^g\Big(1-\a_{X,\F_q;\ell}^{(q_{{\bf n}_m})}x\Big)\Big(1-\b_{X,\F_q;\ell}^{(q_{{\bf n}_m})}x\Big)}{(1-x)(1-q_{{\bf n}_m} x)}
\ea
$$
Therefore, by clearing the common denominator and comparing coefficients of $x^n$, we obtain the following
\begin{thm}\label{thm3.2} With the same notation as above, for $k\geq 0$,
$$\ba
q_{{\bf n}_m}^kb_{X,\F_q;k}^{({\bf n}_m)}-(q_{{\bf n}_m}-1)&q_{{\bf n}_m}^{k-1}b_{X,\F_q;k-1}^{({\bf n}_m)}+q_{{\bf n}_m}^{k-2}b_{X,\F_q;k-2}^{({\bf n}_m)}
=\sum_{\substack{0\leq\ell,0\leq\ell'\leq 2g\\ \ell+\ell'=k}}A_{X,\F_q;\ell'}^{(q_{{\bf n}_m})} \cdot b_{X,\F_q;\ell'}^{({\bf n}_m)}\\
=&b_{X,\F_q;k}^{({\bf n}_m)}+A_{X,\F_q;1}^{(q_{{\bf n}_m})} \cdot b_{X,\F_q;k-1}^{({\bf n}_m)}+\cdots+A_{X,\F_q;2g}^{(q_{{\bf n}_m})} \cdot b_{X,\F_q;k-2g}^{({\bf n}_m)}
\ea$$
where we have set, for $k=-2g,\ldots,-2-1,0$
$$
b_{X,\F_q;k}^{({\bf n}_m)}=\bc 1&k=1\\
0&k=-1,-2,\ldots,-2g.
\ec$$
\end{thm}
\bp Indeed, by clearing up the common denominator on both sides of the relation immediately before this theorem, we get
$$
{(1-x)(1-q_{{\bf n}_m} x)}\cdot B_{X,\F_q}^{({\bf n}_m)}(q_{{\bf n}_m}x)=B_{X,\F_q}^{({\bf n}_m)}(x)
\cdot {\prod_{\ell=1}^g\Big(1-\a_{X,\F_q;\ell}^{(q_{{\bf n}_m})}x\Big)\Big(1-\b_{X,\F_q;\ell}^{(q_{{\bf n}_m})}x\Big)}
$$
That is to say
$$
{(1-x)(1-q_{{\bf n}_m} x)}\sum_{k=0}^\infty b_{X,\F_q;k}^{({\bf n}_m)}q_{{\bf n}_m}^kx^k=\sum_{k=0}^\infty b_{X,\F_q;k}^{({\bf n}_m)}x^k{\prod_{\ell=1}^g\Big(1-\a_{X,\F_q;\ell}^{(q_{{\bf n}_m})}x\Big)\Big(1-\b_{X,\F_q;\ell}^{(q_{{\bf n}_m})}x\Big)}.
$$
By Corollary\,\ref{cor2.5}, let now
$$
\ba
&\frac{P_{X,\F_q}^{(q_{{\bf n}_m})}(x)}{ \a_{X,\F_q}^{(q_{{\bf n}_m})}(0)}=\sum_{\ell=0}^{2g}A_{X,\F_q;\ell}^{(q_{{\bf n}_m})}x^\ell
=\prod_{\ell=1}^g\Big(1-\a_{X,\F_q;\ell}^{(q_{{\bf n}_m})}x\Big)\Big(1-\b_{X,\F_q;\ell}^{(q_{{\bf n}_m})}x\Big)\\
=&\Biggl(\Bigg(\sum_{m=0}^{g-2}\a_{X,\F_q;m}^{(q_{{\bf n}_m})'}\Big(x^{m}+q_{{\bf n}_m}^{(g-1)-m}x^{2(g-1)-m}\Big)+\a_{X,\F_q;g-1}^{(q_{{\bf n}_m})'}x^{g-1}\Bigg)(1-x)(1-q_{{\bf n}_m}x)+(q_{{\bf n}_m}-1)\b_{X,\F_q}^{({\bf n}_m)'}x^g\Biggr)\cdot\\
\ea
$$ 
where 
$$
A_{X,\F_q;\ell}^{(q_{{\bf n}_m})}=(-1)^\ell\sum_{1\leq l_1<\ldots<l_k\leq 2g}\a_{X,\F_q;\ell_1}^{(q_{{\bf n}_m})}\a_{X,\F_q;\ell_2}^{(q_{{\bf n}_m})}\a_{X,\F_q;\ell_k}^{(q_{{\bf n}_m})}
$$
which are given by
$$
A_{X,\F_q;\ell}^{(q_{{\bf n}_m})}=\bc\a_{X,\F_q;0}^{(q_{{\bf n}_m})'}=
1&k=0\\
\a_{X,\F_q;1}^{(q_{{\bf n}_m})'}&k=1\\
\a_{X,\F_q;\ell}^{(q_{{\bf n}_m})'}-(1+q_{{\bf n}_m})\a_{X,\F_q;\ell-1}^{(q_{{\bf n}_m})'}+q_{{\bf n}_m}\a_{X,\F_q;\ell-2}^{(q_{{\bf n}_m})'}&2\leq k\leq g-1\\
(q_{{\bf n}_m}-1)\b_{X,\F_q}^{(q_{{\bf n}_m})'}-(1+q_{{\bf n}_m})\a_{X,\F_q;g-1}^{(q_{{\bf n}_m})'}+q_{{\bf n}_m}\a_{X,\F_q;g-2}^{(q_{{\bf n}_m})'}&k=g\\
q_{{\bf n}_m}^2\a_{X,\F_q;g-3}^{(q_{{\bf n}_m})'}+q_{{\bf n}_m}\a_{X,\F_q;g-1}^{(q_{{\bf n}_m})'}&k=g+1\\
q_{{\bf n}_m}^{k-g+1}\a_{X,\F_q;2g-2-k}^{(q_{{\bf n}_m})'}-(1+q_{{\bf n}_m})q_{{\bf n}_m}^{k-g}\a_{X,\F_q;2g-2-k+1}^{(q_{{\bf n}_m})'}\\
\qquad\qquad +q_{{\bf n}_m}q_{{\bf n}_m}^{k-g-1}\a_{X,\F_q;2g-2-k+2}^{(q_{{\bf n}_m})'}&g+2\leq k\leq 2g-1\\
q_{{\bf n}_m}^g&k=2g
\ec
$$
Then
$$
\Big(1-(q_{{\bf n}_m}+1) x+q_{{\bf n}_m}x^2\Big)\sum_{k=0}^\infty b_{X,\F_q;k}^{({\bf n}_m)}q_{{\bf n}_m}^kx^k=\sum_{\ell=0}^{2g}A_{X,\F_q;\ell}^{(q_{{\bf n}_m})}x^\ell\cdot \sum_{k=0}^\infty b_{X,\F_q;k}^{({\bf n}_m)}x^k.
$$
Or equivalently,
$$\ba
&\sum_{k=0}^\infty b_{X,\F_q;k}^{({\bf n}_m)}q_{{\bf n}_m}^kx^k-(q_{{\bf n}_m}+1) \sum_{k=0}^\infty b_{X,\F_q;k}^{({\bf n}_m)}q_{{\bf n}_m}^kx^{k+1}+q_{{\bf n}_m} \sum_{k=0}^\infty b_{X,\F_q;k}^{({\bf n}_m)}q_{{\bf n}_m}^kx^{k+2}\\
&\hskip 5.60cm=\sum_{k=0}^\infty\sum_{\ell=0}^{2g}A_{X,\F_q;\ell}^{(q_{{\bf n}_m})} \cdot b_{X,\F_q;k}^{({\bf n}_m)}x^{k+\ell}.
\ea$$
Therefore,
$$\ba
q_{{\bf n}_m}^kb_{X,\F_q;k}^{({\bf n}_m)}-(q_{{\bf n}_m}-1)&q_{{\bf n}_m}^{k-1}b_{X,\F_q;k-1}^{({\bf n}_m)}+q_{{\bf n}_m}^{k-2}b_{X,\F_q;k-2}^{({\bf n}_m)}=\sum_{\substack{0\leq\ell,0\leq\ell'\leq 2g\\ \ell+\ell'=k}}A_{X,\F_q;\ell'}^{(q_{{\bf n}_m})} \cdot b_{X,\F_q;\ell'}^{({\bf n}_m)}\\
=&b_{X,\F_q;k}^{({\bf n}_m)}+A_{X,\F_q;1}^{(q_{{\bf n}_m})} \cdot b_{X,\F_q;k-1}^{({\bf n}_m)}+\cdots+A_{X,\F_q;2g}^{(q_{{\bf n}_m})} \cdot b_{X,\F_q;k-2g}^{({\bf n}_m)}
\ea$$
as wanted.
\ep
Furthermore, from the same last relation immediately before this theorem, by replacing $x$ with $x/q_{{\bf n}_m}^k$, we get
$$
\frac{B_{X,\F_q}^{({\bf n}_m)}(q_{{\bf n}_m}x/q_{{\bf n}_m}^k)}{B_{X,\F_q}^{({\bf n}_m)}(x/q_{{\bf n}_m}^k)}=\prod_{\ell=1}^g\frac{\Big(1-\a_{X,\F_q;\ell}^{(q_{{\bf n}_m})}x/q_{{\bf n}_m}^k\Big)\Big(1-\b_{X,\F_q;\ell}^{(q_{{\bf n}_m})}x/q_{{\bf n}_m}^k\Big)}{\Big(1-x/q_{{\bf n}_m}^k\Big)\Big(1-q_{{\bf n}_m} x/q_{{\bf n}_m}^k\Big)}
$$

Taking the product for $k=1,2,\ldots,\infty$, we conclude then
\be\label{eq70}
B_{X,\F_q}^{({\bf n}_m)}(x)=\prod_{k=1}^\infty\prod_{\ell=1}^g\frac{\Big(1-\a_{X,\F_q;\ell}^{(q_{{\bf n}_m})}q_{{\bf n}_m}^{-k}x\Big)\Big(1-\b_{X,\F_q;\ell}^{(q_{{\bf n}_m})}q_{{\bf n}_m}^{-k}x\Big)}{\Big(1-q_{{\bf n}_m}^{-k}x\Big)\Big(1-q_{{\bf n}_m}^{1-k} x\Big)}
\ee
In particular, inserting $x$ with $q_{{\bf n}_{m-1}}^{-n_ms}$, we get
$$
B_{X,\F_q}^{({\bf n}_m)}(T_{{\bf n}_m})=\prod_{k=1}^\infty\frac{\wh\zeta_{X,\F_q}^{\,({\bf n}_m)}(k+s)}{\a_{X,\F_q}^{\,({\bf n}_m)}(0)}.
$$

To use this formula, following\cite{WZ1}, we recall that, for a fixed positive integer $n$,  the \lq\lq $q$-Pochhammer symbol'' $(x,q)_n$ is defined for $x,q\in \C$ as
\be
(x,q)_n=\prod_{m=0}^{n-1}(1-q^mx)
\ee
This can be extended to
\be
(x,q)_\infty=\prod_{m=0}^\infty(1-q^mx)\qquad (\forall |q|<1)
\ee
In our later use, $q$ is supposed to be  $q_{{\bf n}_m}$, which is bigger than 1. Hence we will replace it by its inverse. This then leads to a version of the so-called \lq\lq quantum dilogarithm identity''
\be
\sum_{k=1}^\infty\frac{x^k}{k(q_{{\bf n}_m}^k-1)}=\sum_{k,r\geq 1}\frac{q_{{\bf n}_m}^{-rk}x^k}{k}=\log\frac{1}{\Big(q_{{\bf n}_m}^{-1}x;q_{{\bf n}_m}^{-1}\Big)_\infty}\qquad(|q_{{\bf n}_m}|>1).
\ee
Consequently, we have
\be
B_{X,\F_q}^{({\bf n}_m)}(x)=\frac{\prod_{\ell=1}^g\Big(q_{{\bf n}_m}^{-1}\a_{X,\F_q;\ell}^{(q_{{\bf n}_m})}x;q_{{\bf n}_m}^{-1}\Big)_\infty\Big(q_{{\bf n}_m}^{-1}\b_{X,\F_q;\ell}^{(q_{{\bf n}_m})}x;q_{{\bf n}_m}^{-1}\Big)_\infty}{\Big(q_{{\bf n}_m}^{-1}x;q_{{\bf n}_m}^{-1}\Big)_\infty\Big(x;q_{{\bf n}_m}^{-1}\Big)_\infty}
\ee
Therefore, by \cite{Z}, particularly, Proposition 2 at p.29, we have proved the following
\begin{thm}\label{thm3.3} The number sequence
$\Big\{ b_{X,\F_q;k}^{({\bf n}_m)}\Big\}_k$ defined in \eqref{eq64} is given by
\be
b_{X,\F_q;k}^{({\bf n}_m)}=\sum_{\substack{k_{0,+},k_{0,-},k_{1+},k_{1-},\ldots k_{g+},k_{g-}\geq0\\
k_{0+}+k_{0-}+k_{1+}+k_{1-}+\cdots+ k_{g+}+k_{g-}=k}}
\frac{(-1)^{k_{0+}+k_{0-}}q_{{\bf n}_m}^{\binom{k_{0+}+1}{2}+\binom{k_{0-}}{2}}}{\Big(q_{{\bf n}_m},q_{{\bf n}_m}\Big)_{k_{0+}}\Big(q_{{\bf n}_m},q_{{\bf n}_m}\Big)_{k_{0-}}}\prod_{\ell=1}^g\frac{\Big(\a_{X,\F_q;\ell}^{(q_{{\bf n}_m})}\Big)^{k_{\ell+}}\Big(\b_{X,\F_q;\ell}^{(q_{{\bf n}_m})}\Big)^{k_{\ell-}}}
{\Big(q_{{\bf n}_m},q_{{\bf n}_m}\Big)_{k_{\ell+}}\Big(q_{{\bf n}_m},q_{{\bf n}_m}\Big)_{k_{\ell-}}}\ee
\end{thm}

In the next subsection, we will expose that, when $X$ is an elliptic curve $E$, the invariants $b_{E,\F_q;k}^{({\bf n}_m)}$ are closed related to the beta invariants
$\b_{E,\F_q;k}^{({\bf n}_{m})}$.

\subsection{Relation between $b_{E,\F_q;k}^{({\bf n}_m)}$ and $\b_{E,\F_q;k}^{({\bf n}_{m})}$ for an elliptic curve $E$}

In this subsection, using the techniques developed in \cite{WZ1}, we will give an intrinsic relation between $b_{E,\F_q;k}^{({\bf n}_m)}$ and
$\b_{E,\F_q;k}^{({\bf n}_{m})}$ when $X$ restricts to an elliptic curve $E$ over $\F_q$. The discussion here will be used in the final subsection to prove the ${\bf n}_m$-derived Riemann hypothesis for elliptic curves over finite fields.

Recall that, for general $X/\F_q$,
$$
\b_{X,\F_q}^{({\bf n}_m)}:=\sum_{\substack{k_1,k_2,\ldots, k_p\geq 1\\ k_1+k_2+\cdots+ k_p=n_m}}\frac{(-1)^{p-1}\wh v_{X,\F_q;k_1}^{\,({\bf n}_{m-1})}\wh v_{X,\F_q;k_2}^{\,({\bf n}_{m-1})}\cdots \wh v_{X,\F_q;k_p}^{\,({\bf n}_{m-1})}}{\Big(q_{{\bf n}_{m-1}}^{k_1+k_2}-1\Big)\Big(q_{{\bf n}_{m-1}}^{k_2+k_3}-1\Big)\cdots \Big(q_{{\bf n}_{m-1}}^{k_{p-1}+k_p}-1\Big)}
$$
Here we have replaced and will replace $\wh\zeta_{X,\F_q}^{({\bf n}_m)}(s)$ with $\frac{1}{\a_{X,\F_q}^{({\bf n}_m)}(0)}\wh\zeta_{X,\F_q}^{({\bf n}_m)}(s)$ so that the constant term of the polynomial
$P_{X,\F_q}^{({\bf n}_m)}(T_{{\bf n}_m})$ becomes one. By the comment at the end of \S\ref{sec1.3}, we see that this modification with constant factors is compatible with our discussion below since the closed formula in Theorem\,\ref{thm2.6} is compatible with such a replacement.

According to  $k_p=1,\ldots, n_m$, rewrite $\b_{X,\F_q}^{({\bf n}_m)}$ as a sum
$$
\b_{X,\F_q}^{({\bf n}_m)}=\sum_{\ell=1}^{n_m}\b_{X,\F_q;\ell}^{({\bf n}_m)}\qquad (n_m\geq 1)
$$
where
$$\ba
\b_{X,\F_q;\ell}^{({\bf n}_m)}:=&\sum_{\substack{k_1,k_2,\ldots, k_{p-1}\geq 1, k_p=\ell\\
k_1+k_2+\cdots+ k_p=n_m}}\frac{(-1)^{p-1}\wh v_{X,\F_q;k_1}^{\,({\bf n}_{m-1})}\wh v_{X,\F_q;k_2}^{\,({\bf n}_{m-1})}\cdots \wh v_{X,\F_q;k_p}^{\,({\bf n}_{m-1})}}{\Big(q_{{\bf n}_{m-1}}^{k_1+k_2}-1\Big)\Big(q_{{\bf n}_{m-1}}^{k_2+k_3}-1\Big)\cdots \Big(q_{{\bf n}_{m-1}}^{k_{p-1}+k_p}-1\Big)}\\
=&\wh v_{X,\F_q;\ell}^{\,({\bf n}_{m-1})}\bc 1&\ell=n_m\\
-\sum_{k=1}^{n_m-\ell}\frac{\b_{X,\F_q;k}^{({\bf n}_m-\ell)}}{q_{{\bf n}_{m-1}}^{\ell+k}-1}&1\leq \ell\leq n_m-1\ec
\ea
$$
where for simplicity, we have set
$$
{\bf n}_m-\ell=(n_0,n_2,\ldots,n_m-\ell)
$$
This recursion in turn defines all of the numbers 
$\b_{X,\F_q;\ell}^{({\bf n}_m)}$ (and hence also all of the numbers $\b_{X,\F_q}^{({\bf n}_m)}$). Multiplying this formula by $x^{n_m}$ and summing over all $n_m\geq 0$, we find that the generating functions
$$
b^{({\bf n}_{m-1})}_{X,\F_q;\ell}(x)=\sum_{\substack{n_m=0\\ {\bf n}_m=({\bf n}_{m-1},\,n_m)}}^\infty\b_{X,\F_q;\ell}^{({\bf n}_{m})}x^{n_m}
$$
of the $\b_{X,\F_q;\ell}^{({\bf n}_{m-1})}$ (observe that the sum here actually starts at $n_m=\ell$, so that
$$
b^{({\bf n}_{m-1})}_{X,\F_q;\ell}(x)=O(x^\ell)\qqan b^{({\bf n}_m)}_{X,\F_q;0}(x):\equiv 1)
$$ 
satisfy the identity
\be\label{eq82}
b^{({\bf n}_{m-1})}_{X,\F_q;\ell}(x)=\wh v_{X,\F_q;\ell}^{\,({\bf n}_{m-1})}x^\ell\left(1-\sum_{k=1}^\infty\frac{b_{X,\F_q;k}^{({\bf n}_{m-1})}(x)}{q_{{\bf n}_{m-1}}^{\ell+k}-1}\right)
\qquad(\ell\geq 1)\ee

From now on till the end of this subsection, we restrict ourselves to the case that  $X$ is an  elliptic curve $E$ over $\F_q$. Since $g=1$,  the factor $q_{{\bf n}_m}^{g-1}=1$. So we may omit $\wh~\,$ in all the discussions. For example,
\be\label{eq83}
\wh\zeta_{E,\F_q}^{\,({\bf n}_{m})}(s)=\zeta_{E,\F_q}^{\,({\bf n}_{m})}(s)=\a_{E,\F_q}^{({\bf n}_{m})}(0)\cdot\frac{1-a_{E,\F_q}^{({\bf n}_{m})}q_{{\bf n}_{m-1}}^{-s}-q_{{\bf n}_{m}}q_{{\bf n}_{m-1}}^{-2s}}{(1-q_{{\bf n}_{m}}^{-s})(1-q_{{\bf n}_{m}}^{1-s})}
\ee
where 
$$
a_{E,\F_q}^{({\bf n}_{m})}:=a_{E,\F_q;1}^{({\bf n}_{m})}=\a_{E,\F_q;1}^{({\bf n}_{m})}+\b_{E,\F_q;1}^{({\bf n}_{m})}
$$

With this, particular, by omitting the constant factor $\a_{E,\F_q}^{({\bf n}_{m})}(0)$ at all levels,
we are now ready to state the following
\begin{thm}\label{thm3.4}
Let $E$ be an integral regular elliptic curve over $\F_q$. Define the ${\bf n}_m$-derived Derichlet series
$$
\frZ_{E,\F_q}^{({\bf n}_{m-1})}(s):=\sum_{\substack{n=0\\ {\bf n}_m=({\bf n}_{m-1},n)}}^\infty
\b_{E,\F_q}^{({\bf n}_{m})}q_{{\bf n}_{m-1}}^{-ns}.
$$
Then we have
$$
\frZ_{E,\F_q}^{({\bf n}_{m-1})}(s)=\prod_{k=1}^\infty\zeta_{E,\F_q}^{({\bf n}_{m-1})}(s+k).
$$
Or equivalently,
\be\label{eq86}
\b_{E,\F_q}^{({\bf n}_{m})}=\b_{E,\F_q;n}^{({\bf n}_{m-1})}=b_{E,\F_q;n}^{({\bf n}_{m-1})}=:b_{E,\F_q}^{({\bf n}_{m})}\qquad\forall {\bf n}_{m}=({\bf n}_{m-1},n)
\ee
\end{thm}

\bp
Obviously, both side of \eqref{eq86} are polynomials of $a_{E,\F_q}^{({\bf n}_{m-1})}$. Hence it suffices to verify \eqref{eq86} for  infinitely many special values of $a_{E,\F_q}^{({\bf n}_{m-1})}$ for  fixed ${\bf n}_{m}$ or the same for fixed ${\bf n}_{m-1}$ and $n$. Motivated by \cite{WZ1}, we take
\be\label{eq87}
a_{E,\F_q}^{({\bf n}_{m-1})}:=a_k=q^{k+1}+q^{-k}\qquad(k\in\Z_{\geq 0})
\ee
Accordingly denote by $\b_{E,\F_q}^{({\bf n}_{m}),k}$ and $\b_{E,\F_q;\ell}^{({\bf n}_{m-1}),k}$ the specialization of $\b_{E,\F_q}^{({\bf n}_{m})}$ and $\b_{E,\F_q;\ell}^{({\bf n}_{m-1})}$ to this value of $a_{E,\F_q}^{({\bf n}_{m-1})}:=a_k$ and by $b^{({\bf n}_{m}),k}_{X,\F_q}(x)$ and 
$b^{({\bf n}_{m-1}),k}_{X,\F_q;\ell}(x)$ the corresponding generating series. Then \eqref{eq82}
specializes to the identity
\be\label{eq91}
b^{({\bf n}_{m-1}),k}_{X,\F_q;\ell}(x)=\wh v_{X,\F_q;\ell}^{\,({\bf n}_{m-1}),k}x^\ell\left(1-\sum_{p=1}^\infty\frac{b_{X,\F_q;p}^{({\bf n}_{m-1}),k}(x)}{q_{{\bf n}_{m-1}}^{\ell+p}-1}\right)
\qquad(\ell\geq 1)\ee
where $v_{X,\F_q;\ell}^{\,({\bf n}_{m-1}),k}$ denotes the specialization of $v_{X,\F_q;\ell}^{\,({\bf n}_{m-1})}$ to $a_{E,\F_q}^{({\bf n}_{m-1})}:=a_k$ using \eqref{eq83}. In this way, if we can guess some other of numbers $\wt  \b_{E,\F_q;\ell}^{({\bf n}_{m}),k}$, say $b_{E,\F_q;\ell}^{({\bf n}_{m}),k}$, whose generating function satisfy the same identity, then we automatically have

$$
\b_{E,\F_q;\ell}^{({\bf n}_{m}),k}=\wt \b_{E,\F_q;\ell}^{({\bf n}_{m}),k}=b_{E,\F_q;\ell}^{({\bf n}_{m}),k}.
$$

The reason for looking at the special values \eqref{eq87} is that the relation \eqref{eq70} for this value of $a_{E,\F_q}^{({\bf n}_{m-1})}:=a_k$ implies that the generating function $B_{X,\F_q}^{({\bf n}_{m-1}),k} \ (k\geq 0)$ is given by
\be\label{eq93}
B_{E,\F_q}^{({\bf n}_{m-1}),k}(x)=\prod_{r=1}^\infty\frac{(1-q_{{\bf n}_{m-1}}^{-k-r}x)(1-q_{{\bf n}_{m-1}}^{k+1-r}x)}{(1-q_{{\bf n}_{m-1}}^{-r}x)(1-q_{{\bf n}_{m-1}}^{1-r} x)}=\prod_{j=1}^\infty\frac{1-q_{{\bf n}_{m-1}}^{j}x}{1-q_{{\bf n}_{m-1}}^{-j}x}.
\ee
(in particular, it is a rational function of $x$) and also that the numbers $v_{E,\F_q;\ell}^{\,({\bf n}_{m-1}),k}$ are given by
$$
v_{E,\F_q;\ell}^{\,({\bf n}_{m-1}),k}=\bc(-1)^{\ell+1}q_{{\bf n}_{m-1}}^{\binom{\ell}{2}-k\ell}
\frac{(q_{{\bf n}_{m-1}})_{\ell+k}}{(q_{{\bf n}_{m-1}})_\ell(q_{{\bf n}_{m-1}})_{\ell-1}(q_{{\bf n}_{m-1}})_{k-\ell}}&1\leq\ell\leq k\\
0&\ell>k\ec
$$
as one checks directly using Proposition 2 at p.29 of \cite{Z}. Here and the rest of the section, we use the notation $(x)_{({\bf n}_{m-1},n)}$ for the $q_{{\bf n}_{m-1}}$-Pochhammer symbol 
$$
(x)_{({\bf n}_{m-1},n)}:=\bc 1&n=0\\
(1-x)(1-q_{{\bf n}_{m-1}}x)\cdots(1-q_{{\bf n}_{m-1}}^{n-1}x)&n\geq 1.
\ec
$$
Note that $(x)_{({\bf n}_{m-1},n)}$ is indeed the $q_{{\bf n}_{m-1}}$-Pochhammer symbol, obtained from the $q$-Pochhammer symbol by a simple
replacement of  $q$ with $q_{{\bf n}_{m-1}}$. Hence the $q_{{\bf n}_{m-1}}$-Pochhammer symbol$(x)_{({\bf n}_{m-1},n)}$ satisfies all standard properties for the $q$-Pochhammer symbols. For example,
if we let
$$
\left[{k\atop r}\right]_{({\bf n}_{m-1})}:=\frac{(q_{{\bf n}_{m-1}})_{({\bf n}_{m-1},k)}}{(q_{{\bf n}_{m-1}})_{({\bf n}_{m-1},r)}(q_{{\bf n}_{m-1}})_{{({\bf n}_{m-1},k-r)}}}
$$
be a generalized $q$ or better $q_{{\bf n}_{m-1}}$-binomial coefficient. Naturally, these $q_{{\bf n}_{m-1}}$-binomial coefficients satisify the following  two well-known $q_{{\bf n}_{m-1}}$-versions of the binomial theorem
\be\label{eq97}
\sum_{r=0}^k(-1)^rq_{{\bf n}_{m-1}}^{\binom{r}{2}}\left[{k\atop r}\right]_{({\bf n}_{m-1})}x^r=(x)_{({\bf n}_{m-1},k)}\qqan
\sum_{r=0}^\infty \left[{k+r-1\atop r}\right]_{({\bf n}_{m-1})}x^r=\frac{1}{(x)_{({\bf n}_{m-1},k)}},
\ee
where $k$ denotes an integer $\geq0$. 

\begin{prop}\label{prop3.5} For $k\geq 0$ and $\ell\geq 1$, the generating function $ B_{X,\F_q;\ell}^{({\bf n}_{m-1}),k}(x)$ is a rational function of $x$, equal to 0 if $\ell>k$ and otherwise given by
\be\label{eq98}
B_{X,\F_q;\ell}^{({\bf n}_{m-1}),k}(x)=(-1)^{\ell-1}\frac{(q_{{\bf n}_{m-1}})_{({\bf n}_{m-1},\ell+k)}}{(q_{{\bf n}_{m-1}})_{({\bf n}_{m-1},k)}(q_{{\bf n}_{m-1}})_{({\bf n}_{m-1},\ell-1)}}\frac{\ x^\ell Y_k^{({\bf n}_{m-1}),\ell}(x)\ }{D_k^{({\bf n}_{m-1})}(x)}
\ee
where $D_k^{({\bf n}_{m-1})}(x)$ is defined by the product expansion
$$
D_k^{({\bf n}_{m-1})}(x)=\prod_{j=1}^k(q_{{\bf n}_{m-1}}^j-x)
$$
and where $Y_k^{({\bf n}_{m-1}),\ell}(x)$ is the polynomial of degree $k-\ell$ defined by 
\be\label{eq100}\ba
Y_k^{({\bf n}_{m-1}),\ell}(x):=&\sum_{r=0}^{k-\ell}q_{{\bf n}_{m-1}}^{\binom{r+1}{2}+\binom{k-\ell-r+1}{2}} \left[{k\atop r}\right]_{({\bf n}_{m-1})} \left[{k\atop k-\ell-r}\right]_{({\bf n}_{m-1})}x^r\\
=&{\rm Coefficient\ of\ } Z^{k-\ell}\ {\rm in}\ \prod_{j=1}^k\Big(1+
q_{{\bf n}_{m-1}}^jZ\Big)\Big(1+q_{{\bf n}_{m-1}}^jZx\Big)\ea
\ee
\end{prop}
\bp First note that, easily, the equation of the two expressions in the last relation for $Y_k^{({\bf n}_{m-1}),\ell}(x)$ follows from the first formula of \eqref{eq97}.

To prove this proposition, we take \eqref{eq98}, with $Y_k^{({\bf n}_{m-1}),\ell}(x)$ defined as 0 if $\ell>k$, as the definition of the power series $B_{X,\F_q;\ell}^{({\bf n}_{m-1}),k}(x)$
for all $\ell\geq 1$ and $k\geq 0$ and then prove that these power series satisfy the identity \eqref{eq91}. Inserting equations \eqref{eq93}, \eqref{eq98} and the defining formula \eqref{eq100} for $Y_k^{({\bf n}_{m-1}),\ell}(x)$ into \eqref{eq91}, we conclude that after multiplying both sides by a common factor that the identity to be proved is
\be\label{eq101}\ba
q_{{\bf n}_{m-1}}^{k\ell-\binom{\ell}{2}}Y_k^{({\bf n}_{m-1}),\ell}(x)=&\left[k \atop \ell\right]_{({\bf n}_{m-1})}D_k^{({\bf n}_{m-1})}(x)\\
&+\frac{(q_{{\bf n}_{m-1}})_{({\bf n}_{m-1},k+1)}}{(q_{{\bf n}_{m-1}})_{({\bf n}_{m-1},\ell)}(q_{{\bf n}_{m-1}})_{({\bf n}_{m-1},k-\ell)}}\sum_{p=1}^k\left[k+p\atop k+1\right]_{({\bf n}_{m-1})}Y_k^{({\bf n}_{m-1}),p}(x)\frac{(-x)^p}{q_{{\bf n}_{m-1}}^{\ell+p}-1}.
\ea
\ee
But by  a similar yet simpler partial fractions decomposition  argument as the one we just used above we conclude that $Y_k^{({\bf n}_{m-1}),\ell}(x)$ equals $(q_{{\bf n}_{m-1}}^k+1x)^{-\ell}$ times the coefficient of $Z^{k-\ell}$ in the same product $\prod_{j=1}^k\Big(1-q_{{\bf n}_{m-1}}^jZ\Big)\Big(1-q_{{\bf n}_{m-1}}^jxZ\Big)$ as the one used in the original definition \eqref{eq100}  of $Y_k^{({\bf n}_{m-1}),\ell}(x)$, so that the left hand side of \eqref{eq101} can be written, using the first equation of \eqref{eq97}, as
$$
q_{{\bf n}_{m-1}}^{k\ell-\binom{\ell}{2}}Y_k^{({\bf n}_{m-1}),\ell}(x)={\rm Coefficient\ of\ } Z^{k}\ {\rm in}\ \prod_{j=1}^k\Big(1+
q_{{\bf n}_{m-1}}^jZ\Big)\sum_{s=0}^{k-\ell}q_{{\bf n}_{m-1}}^{\binom{k+1}{2}+\ell s}\left[k \atop {\ell+s}\right]_{({\bf n}_{m-1})}(xZ)^s.
$$
As such, then, \eqref{eq101}  follows immediately form the following lemma by replacing $x$ by $xZ$, multiplying both sides by $\prod_{j=1}^k\Big(1+q_{{\bf n}_{m-1}}^jZ\Big)$ and comparing the coefficients of $Z^k$ on both sides.

\begin{lem}
For fixed $k\geq 0$ and $\ell\geq 1$, define two power series $\cF_1^{({\bf n}_{m-1})}(x)$ and $\cF_2^{({\bf n}_{m-1})}(x)$ by
$$\ba
\cF_1^{({\bf n}_{m-1})}(x):=&(-q_{{\bf n}_{m-1}}x)_{({\bf n}_{m-1},k)}\sum_{p=1}^\infty\left[k+p\atop k+1\right]_{({\bf n}_{m-1})}\frac{(-x)^p}{q_{{\bf n}_{m-1}}^{\ell+p}-1}\\
\cF_2^{({\bf n}_{m-1})}(x):=&\left[k\atop \ell\right]_{({\bf n}_{m-1})}(1+x)^{-1}-\sum_{s=0}^{k-\ell}q_{{\bf n}_{m-1}}^{\binom{s+1}{2}+\ell s}\left[k\atop \ell+s\right]_{({\bf n}_{m-1})}x^s.
\ea$$
Then
\be\label{eq102}
\cF_2^{({\bf n}_{m-1})}(x)=-\frac{(q_{{\bf n}_{m-1}})_{({\bf n}_{m-1},k+1)}}{(q_{{\bf n}_{m-1}})_{({\bf n}_{m-1},\ell)}(q_{{\bf n}_{m-1}})_{({\bf n}_{m-1},k-\ell)}}\cF_1^{({\bf n}_{m-1})}(x)
\ee
\end{lem}

\bp First, from the second equation of \eqref{eq97}, we conclude easily that the power series $\cF_1^{({\bf n}_{m-1})}(x)$ and $\cF_2^{({\bf n}_{m-1})}(x)$ satisfy the simple functional equation
$$\ba
&\Big(1+q_{{\bf n}_{m-1}}x\Big)q_{{\bf n}_{m-1}}^\ell \cF_2^{({\bf n}_{m-1})}(q_{{\bf n}_{m-1}}x)-\Big(1+q_{{\bf n}_{m-1}}^{k+1}x\Big)\cF_2^{({\bf n}_{m-1})}(x)\\
=&(-q_{{\bf n}_{m-1}}x)_{({\bf n}_{m-1},k+1)}\sum_{p=1}^\infty\left[k+p\atop p-1\right]_{({\bf n}_{m-1})}(-x)^p
=\frac{-x}{1+x}.
\ea
$$
and slightly complicated
functional equation, by using telescoping series,
$$\ba
&\Big(1+q_{{\bf n}_{m-1}}x\Big)q_{{\bf n}_{m-1}}^\ell \cF_2^{({\bf n}_{m-1})}(q_{{\bf n}_{m-1}}x)-\Big(1+q_{{\bf n}_{m-1}}^{k+1}x\Big)\cF_2^{({\bf n}_{m-1})}(x)\\
=&\left[k\atop \ell\right]_{({\bf n}_{m-1})}\left(q_{{\bf n}_{m-1}}^\ell-\frac{1+q_{{\bf n}_{m-1}}^{k+1}x}{1+x}\right)-\sum_{s=0}^{k-\ell}q_{{\bf n}_{m-1}}^{\binom{s+1}{2}+\ell s}\left[k\atop \ell+s\right]_{({\bf n}_{m-1})}\left(\Big(q_{{\bf n}_{m-1}}^{\ell+s}-1\Big)-\Big(q_{{\bf n}_{m-1}}^{k-\ell-s}-1\Big)q_{{\bf n}_{m-1}}^{s+1+\ell} x\right)x^s\\
=&\left[k\atop \ell\right]_{({\bf n}_{m-1})}\left(\Big(1-q_{{\bf n}_{m-1}}^{k+1}\Big)\frac{x}{1+x}-\Big(1-q_{{\bf n}_{m-1}}^{\ell}\Big)\right)+\Big(1-q_{{\bf n}_{m-1}}^k\Big)\sum_{s=0}^{k+\ell}\left[k-1\atop \ell+s-1\right]_{({\bf n}_{m-1})}q_{{\bf n}_{m-1}}^{\binom{s+1}{2}+\ell s}x^s\\
&\qquad-\Big(1-q_{{\bf n}_{m-1}}^k\Big)\sum_{s=0}^{k+\ell-1}\left[k-1\atop \ell+s\right]_{({\bf n}_{m-1})}q_{{\bf n}_{m-1}}^{\binom{s+1}{2}+\ell (s+1)}x^{s+1}\\
=&\frac{(q_{{\bf n}_{m-1}})_{({\bf n}_{m-1},k+1)}}{(q_{{\bf n}_{m-1}})_{({\bf n}_{m-1},\ell)}(q_{{\bf n}_{m-1}})_{({\bf n}_{m-1},k-\ell)}}\frac{x}{1+x}\ea$$
Together these implies \eqref{eq102}, since it is easily seen that a power seroes $\cF(x)$ satisfying
$$
\Big(1+q_{{\bf n}_{m-1}}x\Big)q_{{\bf n}_{m-1}}^\ell\cF(q_{{\bf n}_{m-1}}x)=\Big(1+q_{{\bf n}_{m-1}}^{k+1}x\Big)\cF(x)
$$
for some integers $k\geq 0$ and $\ell\geq 1$ must vanish identically. This proves the lemma.
\ep
With the lemma established, as said earlier, we have completed our proof of the proposition as well.
\ep
We finally verify  \eqref{eq70} using this proposition. For this it suffices to show that the sum over $\ell\geq 1$ of the rational functions \eqref{eq98} coincides with the right hand side of \eqref{eq93}. Combining \eqref{eq98} with  the second part of the relations in \eqref{eq100} for
$Y_k^{({\bf n}_{m-1}),\ell}(x)$ and the second equality \eqref{eq97}, we conclude
$$\ba
\frac{1}{x}\frac{D_k^{({\bf n}_{m-1})}(x)}{1-q_{{\bf n}_{m-1}}^{k+1}}\sum_{\ell=1}^\infty B_{X,\F_q;\ell}^{({\bf n}_{m-1}),k}(x)
=&\sum_{\ell=1}^k(-x)^{\ell-1}\left[{k+\ell \atop k+1}\right]_{({\bf n}_{m-1})}Y_k^{({\bf n}_{m-1}),\ell}(x)\\
=&
{\rm Coefficient\ of\ } Z^{k-1}\ {\rm in}\ \frac{\prod_{j=1}^k\Big(1+
q_{{\bf n}_{m-1}}^jZ\Big)}{\big(1+xZ\big)\Big(1+q_{{\bf n}_{m-1}}^{k+1}xZ\Big)}
\ea
$$
But by comparing poles and residues (after making a partial fractions decomposition), we see that
$$
\frac{\prod_{j=1}^k\Big(1+q_{{\bf n}_{m-1}}^jZ\Big)}{\big(1+xZ\big)\Big(1+q_{{\bf n}_{m-1}}^{k+1}xZ\Big)}
=\frac{\prod_{j=1}^k\Big(1-q_{{\bf n}_{m-1}}^jx^{-1}\Big)}{\Big(1-q_{{\bf n}_{m-1}}^{k+1}\Big)\big(1+xZ\big)}
+\frac{\prod_{j=1}^k\Big(1-q_{{\bf n}_{m-1}}^{j-k-1}x^{-1}\Big)}{\Big(1-q_{{\bf n}_{m-1}}^{-k-1}\Big)\Big(1+q_{{\bf n}_{m-1}}^{k+1}xZ\Big)}+P_{k-2}(Z)
$$
where $P_{k-2}(Z)$ is a polynomial of degree $\leq k-2$ in $Z$. It follows that
$$\ba
\sum_{\ell=1}^\infty B_{X,\F_q;\ell}^{({\bf n}_{m-1}),k}(x)
=&\frac{(-x)^k}{D_k^{({\bf n}_{m-1})}(x)}\left(-\prod_{j=1}^k\Big(1-q_{{\bf n}_{m-1}}^jx^{-1}\Big)+q_{{\bf n}_{m-1}}^{k(k+1)}\prod_{j=1}^k\Big(1-q_{{\bf n}_{m-1}}^{j-k-1}x^{-1}\Big)\right)\\
=&-1+\prod_{j=1}^k\frac{1-q_{{\bf n}_{m-1}}^jx}{1-q_{{\bf n}_{m-1}}^{-j}x}=-1+B_{X,\F_q}^{({\bf n}_{m-1}),k}(x),
\ea
$$
where, in the final equality, we have used the relation \eqref{eq93}. Since $B_k^{(0)}=1$, this completes the proof that
the sum of the function $B_{X,\F_q;\ell}^{({\bf n}_{m-1})}(x)$ defined recursively by \eqref{eq91} coincides with the right hand side of \eqref{eq70} and hence, by what we have been said, completes our proof of the theorem. 
\ep

\section{Zeros of ${\bf n}_m$-Derived Zeta Functions}
\subsection{Riemann hypothesis for $(n_0,n_1,\ldots,n_{m-1}, 2)$-derived zeta functions}
In this subsection we prove the following
\begin{thm}[$(2,2,\ldots,2)$-Derived Riemann Hypothesis] \label{thm4.1} Let $X$ be an integral regular projective curve of genus $g$ over $\F_q$.  For each  fixed $(m+1)$-tuple ${\bf n}_m=(n_0,n_1,\ldots, n_{m-1},2)$ of positive integers, 
the ${\bf n}_m$-derived zeta function $\wh \zeta_{X,\F_q}^{\,({\bf n}_{m})}(T_{{\bf n}_{m}})$ of $X$ over $\F_q$ satisfies the Riemann hypothesis, provided that  the Riemann hypothesis for the ${\bf n}_{m-1}$-derived zeta function $\wh \zeta_{X,\F_q}^{\,({\bf n}_{m-1})}(T_{{\bf n}_{m}})$  holds.\\
 In particular, all the $(2,2,\ldots,2)$-derived zeta function $\wh \zeta_{X,\F_q}^{\,(2,2,\ldots,2)}(s)$ of $X$ over $\F_q$ lie on the line $\Re(s)=\frac{1}{2}$. 
\end{thm}
\bp
We start with $m=0$, then this theorem is proved by Yoshida. For details, please refer to Theorem 2.2 of \cite{WRH} and the related discussions. 

To deal with general situation, from  Definition\,\ref{defn1.1}(2), we have
$$\ba
\wh \zeta_{X,\F_q}^{\,({\bf n}_{m})}(T_{{\bf n}_{m}})
=&
q_{{\bf n}_{m-1}}^{g-1}\sum_{a=1}\sum_{\substack{k_1,\ldots,k_p>0\\ k_1+\ldots+k_p=1}}\frac{\wh v_{X,\F_q,k_1}^{({\bf n}_{m-1})}\ldots\wh v_{X,\F_q,k_p}^{({\bf n}_{m-1})}}{\prod_{j=1}^{p-1}(1-q_{{\bf n}_{m-1}}^{k_j+k_{j+1}})} \frac{1}{(1-q_{{\bf n}_{m-1}}^{2s-1+k_{p}})}\wh\zeta_{X,\F_q}^{\,({\bf n}_{m-1})}(2s-1)\\
&+
q_{{\bf n}_{m-1}}^{g-1}\sum_{a=2}\wh\zeta_{X,\F_q}^{\,({\bf n}_{m-1})}(2s)\sum_{\substack{l_1,\ldots,l_r>0\\ l_1+\ldots+l_r=1}}
 \frac{1}{(1-q_{{\bf n}_{m-1}}^{-2s+1+l_{1}})}\frac{\wh v_{X,\F_q,l_1}^{({\bf n}_{m-1})}\ldots\wh v_{X,\F_q,l_r}^{({\bf n}_{m-1})}}{\prod_{j=1}^{r-1}(1-q_{{\bf n}_{m-1}}^{l_j+l_{j+1}})}\\
\ea$$
Hence
$$
\wh \zeta_{X,\F_q}^{\,({\bf n}_{m})}(T_{{\bf n}_{m}})=q_{{\bf n}_{m-1}}^{g-1}\wh v_{X,\F_q,1}^{({\bf n}_{m-1})}\left(\frac{\wh\zeta_{X,\F_q}^{\,({\bf n}_{m-1})}(2s-1)}{(1-q_{{\bf n}_{m-1}}^{2s})}+
\frac{\wh\zeta_{X,\F_q}^{\,({\bf n}_{m-1})}(2s)}{(1-q_{{\bf n}_{m-1}}^{-2s+2})}\right).
 $$
 By Theorem\,\ref{thm4.1}, we may write
 $$
  \wh\zeta_{X,\F_q}^{\,({\bf n}_{m-1})}(s)=\a_{X,\F_q}^{({\bf n}_{m-1})}(0)\frac{\prod_{\ell=1}^g(1-a_{X,\F_q,{\bf n}_{m-1},\ell}T_{{\bf n}_{m-1}}-q_{{\bf n}_{m-1}}T_{{\bf n}_{m-1}}^2)}{(1-T_{{\bf n}_{m-1}})(1-q_{{\bf n}_{m-1}}T_{{\bf n}_{m-1}})T_{{\bf n}_{m-1}}^{g-1}}
 $$
Furthermore, by our condition on the zeros of $(m-1)$-th derived zeta function $ \wh\zeta_{X,\F_q}^{\,({\bf n}_{m-1})}(s)$, we have, for all $1\leq\ell\leq g$,
  \be
  1-a_{X,\F_q;{\bf n}_{m-1},\ell}T_{{\bf n}_{m-1}}-q_{{\bf n}_{m-1}}T_{{\bf n}_{m-1}}^2)=\Big(1-\a_{X,\F_q;{\bf n}_{m-1},\ell}T_{{\bf n}_{m-1}}\Big)\Big(1-\b_{X,\F_q;{\bf n}_{m-1},\ell}T_{{\bf n}_{m-1}}\Big)
  \ee
  and 
  \be\label{eq40}
  b_{X,\F_q;{\bf n}_{m-1},\ell}=\ov a_{X,\F_q;{\bf n}_{m-1},\ell}\qqan |a_{X,\F_q;{\bf n}_{m-1},\ell}|=\sqrt {q_{{\bf n}_{m-1}}}.
  \ee
  Therefore,
  \be
  \wh \zeta_{X,\F_q}^{\,({\bf n}_{m})}(T_{{\bf n}_{m}})=0\Llra
  \frac{\wh\zeta_{X,\F_q}^{\,({\bf n}_{m-1})}(2s-1)}{(1-q_{{\bf n}_{m-1}}^{2s})}+
\frac{\wh\zeta_{X,\F_q}^{\,({\bf n}_{m-1})}(2s)}{(1-q_{{\bf n}_{m-1}}^{-2s+2})}=0
  \ee
  This is equivalent to,
  \be
   \frac{\wh\zeta_{X,\F_q}^{\,({\bf n}_{m-1})}(2s-1)}{(T_{{\bf n}_{m}}^{-1}-1)}=
\frac{\wh\zeta_{X,\F_q}^{\,({\bf n}_{m-1})}(2s)}{(1-q_{{\bf n}_{m}}T_{{\bf n}_{m}})}.
  \ee
  or better,
  $$\ba
 & (1-q_{{\bf n}_{m}}T_{{\bf n}_{m}})\frac{\prod_{\ell=1}^g\Big(1-\a_{X,\F_q;{\bf n}_{m-1},\ell}q_{{\bf n}_{m-1}}T_{{\bf n}_{m}}\Big)\Big(1-\b_{X,\F_q;{\bf n}_{m-1},\ell}q_{{\bf n}_{m-1}}T_{{\bf n}_{m}}\Big)}{(1-q_{{\bf n}_{m-1}}T_{{\bf n}_{m}})(1-q_{{\bf n}_{m-1}}q_{{\bf n}_{m-1}}T_{{\bf n}_{m}})q_{{\bf n}_{m-1}}^{(g-1)}T_{{\bf n}_{m}}^{g-1}}\\
 =& (T_{{\bf n}_{m}}^{-1}-1)\frac{\prod_{\ell=1}^g\Big(1-\a_{X,\F_q;{\bf n}_{m-1},\ell}T^{-1}_{{\bf n}_{m}}\Big)\Big(1-\b_{X,\F_q;{\bf n}_{m-1},\ell}T^{-1}_{{\bf n}_{m}}\Big)}{(1-T^{-1}_{{\bf n}_{m}})(1-q_{{\bf n}_{m-1}}T^{-1}_{{\bf n}_{m}})T_{{\bf n}_{m}}^{-(g-1)}}\\
  \ea
  $$
  This latest condition is equivalent to
  $$\ba
 & \frac{\prod_{\ell=1}^g\Big(1-\a_{X,\F_q;{\bf n}_{m-1},\ell}q_{{\bf n}_{m-1}}T_{{\bf n}_{m}}\Big)\Big(1-\b_{X,\F_q;{\bf n}_{m-1},\ell}q_{{\bf n}_{m-1}}T_{{\bf n}_{m}}\Big)}{(1-q_{{\bf n}_{m-1}}T_{{\bf n}_{m}})q_{{\bf n}_{m-1}}^{(g-1)}}\\
 =& -\frac{\prod_{\ell=1}^g\Big(T_{{\bf n}_{m}}-\a_{X,\F_q;{\bf n}_{m-1},\ell}\Big)\Big(T_{{\bf n}_{m}}-b_{X,\F_q;{\bf n}_{m-1},\ell}\Big)}{(T_{{\bf n}_{m}}-q_{{\bf n}_{m-1}})T_{{\bf n}_{m}}}\\
  \ea
  $$
    
 To facilitate our ensuing discussion, now we recall the following elementary
  \begin{lem}[Yoshida]\label{lem4.2} Fix a real number $q>1$. Let $\a,\,\b\in \C$ and write $c=\a+\b$. Assume that $\a\b=q$ and that $c\in \R$ satisfies $|c|\leq q+1$. Then for $w\in \C$, we have
$$
|w-\a|\cdot|w-\b|=|1-\a w|\cdot|1-\b w|\ \bc
&\hskip -0.3cm>1\qquad{\rm if}\quad |w|<1\\
&\hskip -0.3cm <1\qquad{\rm if}\quad |w|>1.
\ec
$$
\end{lem}
The interested reader may find a proof of this lemma, together with a generalization  from  Lemma 4.5 of \cite{WRH}. 

Hence, by applying Lemma\,\ref{lem4.2} to $w=q_{{\bf n}_{m-1}}^{-1}T_{{\bf n}_{m}}$, we get

$$\ba
&\left|\frac{\prod_{\ell=1}^g\Big(T_{{\bf n}_{m}}-\a_{X,\F_q;{\bf n}_{m-1},\ell}\Big)\Big(T_{{\bf n}_{m}}-b_{X,\F_q;{\bf n}_{m-1},\ell}\Big)}{(T_{{\bf n}_{m}}-q_{{\bf n}_{m-1}})T_{{\bf n}_{m}}}\right|\\
=
 & \left|\frac{\prod_{\ell=1}^g\Big(1-\a_{X,\F_q;{\bf n}_{m-1},\ell}q_{{\bf n}_{m-1}}T_{{\bf n}_{m}}\Big)\Big(1-\b_{X,\F_q;{\bf n}_{m-1},\ell}q_{{\bf n}_{m-1}}T_{{\bf n}_{m}}\Big)}{(1-q_{{\bf n}_{m-1}}T_{{\bf n}_{m}})q_{{\bf n}_{m-1}}^{(g-1)}}\right|
 \\
 =& \left|\frac{\prod_{\ell=1}^g\Big(q_{{\bf n}_{m-1}}T_{{\bf n}_{m}}-\a_{X,\F_q;{\bf n}_{m-1},\ell}\Big)\Big(q_{{\bf n}_{m-1}}T_{{\bf n}_{m}}-\b_{X,\F_q;{\bf n}_{m-1},\ell}\Big)}{(1-q_{{\bf n}_{m-1}}T_{{\bf n}_{m}})q_{{\bf n}_{m-1}}^{(g-1)}}\right|\\
&\hskip 2.80cm\times \bc<1&\qquad{\rm if}\quad |q_{{\bf n}_{m-1}}T_{{\bf n}_{m}}|<1\\
                >1&\qquad{\rm if}\quad |q_{{\bf n}_{m-1}}T_{{\bf n}_{m}}|>1.
\ec\\
=& \left|\frac{\prod_{\ell=1}^g\Big(q_{{\bf n}_{m-1}}T_{{\bf n}_{m}}\a_{X,\F_q;{\bf n}_{m-1},\ell}^{-1}-1\Big)\Big(q_{{\bf n}_{m-1}}T_{{\bf n}_{m}}\b_{X,\F_q;{\bf n}_{m-1},\ell}^{-1}-1\Big)}{(1-q_{{\bf n}_{m-1}}T_{{\bf n}_{m}})q_{{\bf n}_{m-1}}^{(-1)}}\right|\\
&\hskip 2.80cm\times \bc<1&\qquad{\rm if}\quad |q_{{\bf n}_{m-1}}T_{{\bf n}_{m}}|<1\\
                >1&\qquad{\rm if}\quad |q_{{\bf n}_{m-1}}T_{{\bf n}_{m}}|>1.
\ec\\
=& \left|\frac{\prod_{\ell=1}^g\Big(T_{{\bf n}_{m}}\a_{X,\F_q;{\bf n}_{m-1},\ell}-1\Big)\Big(T_{{\bf n}_{m}}\b_{X,\F_q;{\bf n}_{m-1},\ell}-1\Big)}{(1-q_{{\bf n}_{m-1}}T_{{\bf n}_{m}})q_{{\bf n}_{m-1}}^{(-1)}}\right|\\
&\hskip 2.80cm\times \bc<1&\qquad{\rm if}\quad |q_{{\bf n}_{m-1}}T_{{\bf n}_{m}}|<1\\
                >1&\qquad{\rm if}\quad |q_{{\bf n}_{m-1}}T_{{\bf n}_{m}}|>1.
\ec\\
  \ea
  $$
 This leads to a contradiction, by comparing the first and the last expressions.
  Therefore, the reciprocity roots of the ${\bf n}_m$-derived zeta function $\wh \zeta_{X,\F_q}^{\,({\bf n}_{m})}(T_{{\bf n}_{m}})$ lie neither in the region of neither $|T_{{\bf n}_{m}}q_{{\bf n}_{m-1}}|<1$, nor in the region $|T_{{\bf n}_{m}}q_{{\bf n}_{m-1}}|>1$. Therefore, all reciprocity zeros of  the ${\bf n}_m$-derived zeta function $\wh \zeta_{X,\F_q}^{\,({\bf n}_{m})}(T_{{\bf n}_{m}})$ satisfy the condition
$$
 |\a_{X,\F_q;{\bf n}_{m},\ell}|=\sqrt {q_{{\bf n}_{m}}}\qquad(\forall \ell=1,\ldots,2g).
$$ 
This is exactly what we are wanted.
\ep

\subsection{${\bf n}_m$-derived Riemann hypothesis for elliptic curves over finite fields}
In this subsection, we prove the following
\begin{thm}[${\bf n}_m$-Derived Riemann Hypothesis for Elliptic Curves over Finite Fields]
Let $E$ be an integral regular elliptic curve over $\F_q$.  For each  fixed $(m+1)$-tuple ${\bf n}_m=(n_0,n_1,\ldots,n_m)$ of positive integers, 
the ${\bf n}_m$-derived zeta function $\wh \zeta_{E,\F_q}^{\,({\bf n}_{m})}(T_{{\bf n}_{m}})$ of $E$ over $\F_q$ satisfies the Riemann hypothesis.
\end{thm}
\bp
We use an induction in $m$. When $m=0$, this is proved in \cite{WZ1} by Zagier and myself.

Assume that the ${\bf n}_{m-1}$-derived zeta function $\wh \zeta_{E,\F_q}^{\,({\bf n}_{m-1})}(T_{{\bf n}_{m-1}})$ of $E$ over $\F_q$ satisfies the Riemann hypothesis. We examine the case for
the ${\bf n}_{m}$-derived zeta function $\wh \zeta_{E,\F_q}^{\,({\bf n}_{m})}(T_{{\bf n}_{m}})$, using the techniques developed in \cite{WZ1}.
This now becomes rather easy. Indeed, by directly applying Theorems\,\ref{thm3.2} and \ref{thm3.4} we have the following
\begin{thm}\label{thm4.4} Let $E$ be an integral regular elliptic curve over $\F_q$.  For each  fixed $(m+1)$-tuple ${\bf n}_m=(n_0,n_1,\ldots,n_m)$ of positive integers, the ${\bf n}_m$-derived beta invariants are determined by the following recursion relations:
\be\label{eq112}
\Big(q_{{\bf n}_{m-1}}^{n_m}-1\Big)\,\b_{X,\F_q}^{({\bf n}_m)}=\Big(q_{{\bf n}_{m-1}}^{n_m}+q_{{\bf n}_{m-1}}^{n_m-1}-a_{X,\F_q}^{({\bf n}_{m-1})}\Big)\,\b_{X,\F_q}^{({\bf n}_m)-1}-\Big(q_{{\bf n}_{m-1}}^{n_m-1}-q_{{\bf n}_{m-1}}\Big)\,\b_{X,\F_q}^{({\bf n}_m)-2}
\ee
together with the initial conditions $\b_{X,\F_q}^{({\bf n}_{m-1},0)}=1$ and $\b_{X,\F_q}^{({\bf n}_{m-1},-1)}=0$.
\end{thm}
\begin{cor}\label{cor4.5} We have
\be
1<\frac{\b_{X,\F_q}^{({\bf n}_{m-1},n)}}{\b_{X,\F_q}^{({\bf n}_{m-1},n-1)}}
<\frac{q_{{\bf n}_{m-1}}^{n/2}+1}{q_{{\bf n}_{m-1}}^{n/2}-1}
\ee
\end{cor}
\bp
we prove this corollary using an induction in $n$.

(1) $n=1$. Recall  that, when $n_m=1$, by Example\,\ref{exam1.2}, we have
\be
\wh \zeta_{E,\F_q}^{\,({\bf n}_{m})}(s)=\wh \zeta_{E,\F_q}^{\,({\bf n}_{m-1})}(s).
\ee
This implies that, by the inductive hypothesis in terms of $m$,
\be
\left|a_{E,\F_q}^{\,({\bf n}_{m-1},1)}\right|\leq 2\sqrt {q_{({\bf n}_{m-1},1)}}=2\sqrt {q_{{\bf n}_{m-1}}}.
\ee
Furthermore, by \eqref{eq112},
\be\label{eq112}
\Big(q_{{\bf n}_{m-1}}^{2}-1\Big)\,\b_{X,\F_q}^{({\bf n}_{m-1},2)}=\Big(q_{{\bf n}_{m-1}}^{2}+q_{{\bf n}_{m-1}}-a_{X,\F_q}^{({\bf n}_{m-1})}\Big)\,\b_{X,\F_q}^{({\bf n}_{m-1},1)}-\Big(q_{{\bf n}_{m-1}}-q_{{\bf n}_{m-1}}\Big)\cdot 1
\ee
Hence
$$\ba
1<&
\frac{q_{{\bf n}_{m-1}}^{2}+q_{{\bf n}_{m-1}}-2\sqrt {q_{{\bf n}_{m-1}}}}{q_{{\bf n}_{m-1}}^{2}-1}\\
\leq &\frac{\b_{X,\F_q}^{({\bf n}_{m-1},2)}}{\b_{X,\F_q}^{({\bf n}_{m-1},1)}}
=\frac{q_{{\bf n}_{m-1}}^{2}+q_{{\bf n}_{m-1}}-a_{X,\F_q}^{({\bf n}_{m-1})}}{q_{{\bf n}_{m-1}}^{2}-1}\\
\leq &\frac{q_{{\bf n}_{m-1}}^{2}+q_{{\bf n}_{m-1}}+2\sqrt {q_{{\bf n}_{m-1}}}}{q_{{\bf n}_{m-1}}^{2}-1}
<\frac{q_{{\bf n}_{m-1}}+1}{q_{{\bf n}_{m-1}}-1}
\ea
$$

(2)
Inductively, assume that the assertions hold when $n\leq k-1$, we want to verify the assertion for $n$. Recall that, by Theorem\,\ref{thm4.4}, we have
\be\Big(q_{{\bf n}_{m-1}}^{n}-1\Big)\,\b_{X,\F_q}^{({\bf n}_{m-1},n)}=\Big(q_{{\bf n}_{m-1}}^{n}+q_{{\bf n}_{m-1}}^{n-1}-a_{X,\F_q}^{({\bf n}_{m-1})}\Big)\,\b_{X,\F_q}^{({\bf n}_{m-1},n-1)}-\Big(q_{{\bf n}_{m-1}}^{n-1}-q_{{\bf n}_{m-1}}\Big)\,\b_{X,\F_q}^{({\bf n}_{m-1},n-2)}
\ee
Therefore, by the inductive hypothesis in $n$,
$$\ba
\frac{\b_{X,\F_q}^{({\bf n}_{m-1},n)}}{\ \b_{X,\F_q}^{({\bf n}_{m-1},n-1)}\ }=&\frac{q_{{\bf n}_{m-1}}^{n}+q_{{\bf n}_{m-1}}^{n-1}-a_{X,\F_q}^{({\bf n}_{m-1})}}{q_{{\bf n}_{m-1}}^{n}-1}\,-\frac{q_{{\bf n}_{m-1}}^{n-1}-q_{{\bf n}_{m-1}}}{q_{{\bf n}_{m-1}}^{n}-1}\,\frac{\ \b_{X,\F_q}^{({\bf n}_{m-1},n-2)}\ }{\b_{X,\F_q}^{({\bf n}_{m-1},n-1)}}\\
>&\frac{q_{{\bf n}_{m-1}}^{n}+q_{{\bf n}_{m-1}}^{n-1}-a_{X,\F_q}^{({\bf n}_{m-1})}}{q_{{\bf n}_{m-1}}^{n}-1}\,-\frac{q_{{\bf n}_{m-1}}^{n-1}-q_{{\bf n}_{m-1}}}{q_{{\bf n}_{m-1}}^{n}-1}\\
=&\frac{q_{{\bf n}_{m-1}}^{n}+q_{{\bf n}_{m-1}}-a_{X,\F_q}^{({\bf n}_{m-1})}}{q_{{\bf n}_{m-1}}^{n}-1}>1\\
\ea
$$
and, by the inductive hypothesis in $n$ again,
$$\ba
&\frac{\b_{X,\F_q}^{({\bf n}_{m-1},n)}}{\ \b_{X,\F_q}^{({\bf n}_{m-1},n-1)}\ }-\frac{q_{{\bf n}_{m-1}}^{(n-1)/2}+1}{q_{{\bf n}_{m-1}}^{(n-1)/2}-1}\\
<&\frac{\Big(q_{{\bf n}_{m-1}}^{n}+q_{{\bf n}_{m-1}}^{n-1}-a_{X,\F_q}^{({\bf n}_{m-1})}\Big)-\Big(q_{{\bf n}_{m-1}}^{n/2}+1\Big)^2}{q_{{\bf n}_{m-1}}^{n}-1}
-\frac{q_{{\bf n}_{m-1}}^{n-1}-q_{{\bf n}_{m-1}}}{q_{{\bf n}_{m-1}}^{n}-1}\,\frac{q_{{\bf n}_{m-1}}^{(n-1)/2}-1}{q_{{\bf n}_{m-1}}^{(n-1)/2}+1}\\
<&\frac{1}{q_{{\bf n}_{m-1}}^{n}-1}\left(q_{{\bf n}_{m-1}}^{n-1}+\Big(q_{{\bf n}_{m-1}}+1\Big)-
2q_{{\bf n}_{m-1}}^{n/2}-1-\Big(q_{{\bf n}_{m-1}}^{n-1}-q_{{\bf n}_{m-1}}\Big)\,\frac{q_{{\bf n}_{m-1}}^{(n-1)/2}-1}{q_{{\bf n}_{m-1}}^{(n-1)/2}+1}
\right)\\
=&-\frac{1}{q_{{\bf n}_{m-1}}^{n}-1}\frac{2\Big(q_{{\bf n}_{m-1}}^{n-1}-q_{{\bf n}_{m-1}}^{n/2}\Big)\Big(q_{{\bf n}_{m-1}}^{1/2}-1\Big)}{q_{{\bf n}_{m-1}}^{(n-1)/2}+1}<0
\ea
$$
as wanted.
\ep
Now we are ready to verify the ${\bf n}_m$-derived Riemann hypothesis for elliptic curve $E$ over $\F_q$. Indeed, by Theorem\,\ref{thm2.7}, the special counting miracle and Corollary\,\ref{cor2.5}, modulo the constant factor $\a_{X,\F_q}^{({\bf n}_{m-1})}(0)$  as commented at the end of \S\ref{sec2.5},
we have
\be
Z_{X,\F_q}^{\,({\bf n}_m)}(T_{{\bf n}_m})
=\b_{X,\F_q}^{({\bf n}_{m}-1)}\frac{1-\left(\Big(q_{{\bf n}_m}+1\Big)-\Big(q_{{\bf n}_m}-1\Big)\frac{\b_{X,\F_q}^{({\bf n}_{m-1},n_m)}}{\ \b_{X,\F_q}^{({\bf n}_{m-1},n_m-1)}\ }\right)T_{{\bf n}_m}+q_{{\bf n}_m}T_{{\bf n}_m}^2}{\Big(1-T_{{\bf n}_m}\Big)\Big(1-q_{{\bf n}_m}T_{{\bf n}_m}\Big)}
\ee
By Corollary\,\ref{cor4.5}, the discriminant of the quadratic numerator is negative. Hence, by Lemma\,\ref{lem1.2}, the ${\bf n}_m$-derived Riemann hypothesis holds for elliptic curve $E$ over $\F_q$.
\ep

We end this paper with the following comments. First, the Positivity Conjecture is supported by what we call the asymptotic positivity claiming that when $q$ is sufficiently large,  all ${\bf n}_m$-derived alpha and beta invariants are strictly positive. Secondly,
to establish the general ${\bf n}_m$-derived Riemann hypothesis, it appears that the Positivity Conjecture plays a central role.
Thirdly, the construction of the ${\bf n}_m$-zeta functions admits counterparts for   number fields and  for $L$-functions.  Finally, there is also a question on the relations among various ${\bf n}_m$-derived zeta zeros. All these will be discussed elsewhere in due courses.

 Lin WENG
 
  Graduate School of Mathematics,
Kyushu University, Fukuoka 819-0395, Japan
 
 \& Institute for Fundamental Research,
 $L$-Academy

 weng@math.kyushu-u.ac.jp

\end{document}